\documentclass[12pt]{amsart}
\newtheorem{theorem}{Theorem}[section]
\newtheorem{lemma}[theorem]{Lemma}
\newtheorem{proposition}[theorem]{Proposition}

\begin{document}

\title{Bott-Chern-Aeppli, Dolbeault and Frolicher on Compact Complex 3-folds}
\author[A. McHugh]{Andrew McHugh}
\address{ Department of Computer Science and Mathematics, Penn State-Harrisburg, Middletown, PA \\ {\tt axm964@psu.edu}\  , {\tt andrewmchugh@snet.net}\ }
\begin{abstract}
We give the complete Bott-Chern-Aeppli cohomology for compact complex 3-folds  in terms of Dolbeault, Frolicher, a bi-degree DeRham-like type of cohomology, $K^{p,q}$, defined  as 
$$ K^{p,q}=\frac{ker( \partial ) \cap ker( {\bar{\partial}}) }{im( \partial )\cap ker( {\bar{\partial}} )+im( {\bar{\partial}})\cap ker( \partial )}$$ and ${\check{H}}^1({\mathcal{PH}})$. 
(Here $\mathcal{PH}$ is the sheaf of phuri-harmonic functions.) We then work out the complete Bott-Chern-Aeppli cohomology in some examples.   We give the Bott-Chern-Aeppli cohomology for a hypothetical complex structure on $S^6$  in terms of Dolbeault and Frolicher. We also give the Bott-Chern-Aeppli cohomology on a Calabi-Eckman 3-fold concurring with the calculations of Angella and Tomassini\cite{AngellaAndTomassini}. Finally, we show agreement of our results with the calculation by Angella\cite{Angella} of the Bott-Chern-Aeppli cohomology for small Kuranishi deformations of the Iwasawa manifold.
\end{abstract} 

\maketitle

\section{Introduction}
On a differentiable manifold, the deRham cohomology, 
$$ H^{p} = \frac{ker( d :  
C^{\infty \ p} \rightarrow  C^{\infty \  p+1}  ) }
{im( d :  C^{\infty \   p-1} \rightarrow  C^{\infty \ p})} $$ depends only on the topological structure of the manifold.  (Here, $C^{\infty \ p}$ denotes the {\em complex-valued} p-forms,  as this  is more convenient for the following discussion of cohomology on complex manifolds. )

 On a complex manifold, we also have the Dolbeault cohomology, 
$$ H^{p,q}_{\bar{\partial}} = \frac{ker( {\bar{\partial}} :  
C^{\infty \ p,q} \rightarrow  C^{\infty \  p,q+1}  ) }
{ im( {\bar{\partial}} :  C^{\infty \  p,q-1} \rightarrow  C^{\infty \ p,q}  ) } $$ depends only on the complex structure of the manifold.   When a compact complex manifold is Kahler or satisfies the $\partial {\bar{\partial}}$-lemma, we have the Hodge decomposition, $   H^{r} = \bigoplus_{p+q =r} H^{p,q}_{\bar{\partial}} $.

On a complex manifold of dimension $n$, there are also the Frolicher spectral sequences, $E^{p,q}_r$.  For $r=1$, $E^{p,q}_1  = H^{p,q}_{\bar{\partial}}$ and we have the sequence for $r=1$ and $0\leq q \leq n$,
$$ 0 \rightarrow H^{0,q}_{\bar{\partial}} \stackrel{\partial}{\rightarrow}  H^{1,q}_{\bar{\partial}} \stackrel{\partial}{\rightarrow}  H^{2,q}_{\bar{\partial}} \stackrel{\partial}{\rightarrow} 
\ldots  \stackrel{\partial}{\rightarrow}  H^{n,q}_{\bar{\partial}}  \stackrel{\partial}{\rightarrow} 0  \ \  . $$
So, for $r=2$, we have the definition, $$E^{p,q}_2  = ker(H^{p,q}_{\bar{\partial}} \stackrel{\partial}{\rightarrow} H^{p+1,q}_{\bar{\partial}})  / im(H^{p-1,q}_{\bar{\partial}} \stackrel{\partial}{\rightarrow} H^{p,q}_{\bar{\partial}})  \  . $$
There are further definitions for the Frolicher spectral sequences for $r \geq 2$ and $E^{p,q}_r$ for $r \geq 3$, which we will not be using in this article.

The Aeppli cohomology of a complex manifold is defined by the vector spaces (see Aeppli{\cite{Aeppli}} and also Angella\cite{Angella} and Popovici\cite{Popovici}) :
$$ H^{p,q}_A = \frac{ker( \partial {\bar{\partial}} :  
C^{\infty \ p,q} \rightarrow  C^{\infty \  p+1,q+1}  ) }
{(im(\partial :  C^{\infty \   p-1,q} \rightarrow  C^{\infty \ p,q}) + im( {\bar{\partial}} :  C^{\infty \  p,q-1} \rightarrow  C^{\infty \ p,q}  ) )} $$

The Bott-Chern cohomology of a complex manifold is defined by the vector spaces (see Bott and Chern\cite{BottChern} and also Angella\cite{Angella} and Popovici \cite{Popovici}) :
$$ H^{p,q}_{BC} = \frac{ker( \partial :  
C^{\infty \ p,q} \rightarrow  C^{\infty \  p+1,q}  ) \cap ker( {\bar{\partial}} :  
C^{\infty \ p,q} \rightarrow  C^{\infty \  p,q+1}  ) }
{im(\partial {\bar{\partial}} :  C^{\infty \   p-1,q-1} \rightarrow  C^{\infty \ p,q} )} $$

 On compact complex manifolds, there is a harmonic theory due to Schweitzer{\cite{Schweitzer}} for each of these cohomologies which ensures that they are finite dimensional complex vector spaces.  Schweitzer's harmonic theory shows that the two cohomologies are dual to each other. Let $h^{p,q}_A = dim(H^{p,q}_A)$ and $h^{p,q}_{BC} = dim(H^{p,q}_{BC})$ .  We have then(see Schweitzer\cite{Schweitzer} and also Angella\cite{Angella} and Popovici{\cite{Popovici}}) that $h^{p,q}_A = h^{q,p}_A$, $h^{p,q}_{BC} = h^{q,p}_{BC}$ 
and $h^{p,q}_{A} = h^{n-p,n-q}_{BC}$.  We mention as a historical note that results on compact complex manifolds about finiteness and duality between Aeppli and Bott-Chern cohomology also appear in Bigolin\cite{Bigolin}.

Bott-Chern/Aeppli cohomology has been studied extensively by a number of mathematicians.  
Popovici\cite{Popovici} utilizes Aeppli cohomology, in particular, $H^{n-1,n-1}_A$, to study Gauduchon metrics on complex manifolds.
Tseng and Yau\cite{Tseng-Yau} point out the importance of understanding Bott-Chern/Aeppli cohomology, in particular, $H^{2,2}_{BC}$, for the study of Strominger's system of supersymmetric equations in type IIB theory on complex 3-folds.   

For compact complex surfaces, the deRham and Dolbeault cohomology determine the Bott-Chern and Aeppli cohomology. This is due to a result of Teleman\cite{Teleman}.  (See the discussion in MathOverflow \cite{Daniele}.  Bott-Chern cohomology of compact complex surfaces have also been calculated by Angella, Dloussky and Tomassini\cite{AngellaDlousskyAndTomassini})  
 
Angella\cite{Angella} has given an example of two different complex 3-folds with the same deRham and Dolbeault cohomologies but different Bott-Chern/Aeppli cohomologies.  
In this article we give the complete Bott-Chern-Aeppli cohomology for compact complex 3-folds  in terms of $H^{p,q}_{\bar{\partial}}$, $E^{3,1}_2$, ${\check{H}}^1({\mathcal{PH}})$ and a bi-degree DeRham-like type of cohomology, $K^{p,q}$, defined as
$$ K^{p,q} = \frac{ker( \partial :  
C^{\infty \ p,q} \rightarrow  C^{\infty \  p+1,q}  ) \cap ker( {\bar{\partial}} :  
C^{\infty \ p,q} \rightarrow  C^{\infty \  p,q+1}  ) }
{im( \partial )\cap ker( {\bar{\partial}} ) + im( {\bar{\partial}})\cap ker( \partial )} $$
We set $k^{p,q} = dim_{\mathbf C} (K^{p,q} )$.  
It can easily be shown that  $k^{p,q} \leq b^{p+q}$, $k^{p,q} \leq h_{\bar{\partial}}^{p,q}$, $k^{p,q} \leq h_{\partial}^{p,q} = h_{\bar{\partial}}^{q,p}$.   We also see that $K^{p,q}$ embeds into $H^{p,q}_A$ by the obvious map, $[\kappa]_K \mapsto [\kappa]_A$. This is easily seen to be one-to-one. Thus $k^{p,q} \leq h^{p,q}_A$.  
Since 
$Im(\partial {\bar{\partial}}) \subseteq (Im(\partial) + Im({\bar{\partial}}))  $, 
we have $k^{p,q} \leq h^{p,q}_{BC}$ from the definitions of $H^{p,q}_{BC}$ and $K^{p,q}$. We also easily have from the definition of $K^{p,q}$ that $k^{p,q} = k^{q,p}$.

Our main result is given in the following table:
 \begin{center}
\begin{tabular}{c|c|c|c|c|}
\multicolumn{5}{c}{\bf{Bott-Chern cohomology on a Compact Complex 3-fold}}  \\
\ $h^{p,q}_{BC}$ & $q= 0$ & $ q=1$ & $q=2$ & $ q=3$ \\
\hline 
$p=0$ & \  \  \  \   1 \ \  \ \  &  $k^{1,0}$  & $h^{2,0}_{\bar{\partial}} $ & $k^{3,0}$ \\ 
& \  &  \  &  \   & \   \\ 
$p=1$ & $k^{1,0}$ & \  \ \  \ \ ${\check{h}}^1({\mathcal{PH}})$ \ \  & $ h^{1,2}_{\bar{\partial}}  + h^{3,1}_2 - k^{2,0} $ &  $h^{0,2}_{\bar{\partial}} $ \\ 
& \  &  \  &  \   & \   \\ 
$p=2$ & $h^{2,0}_{\bar{\partial}} $ & $ h^{1,2}_{\bar{\partial}}  + h^{3,1}_2 - k^{2,0} $ &  \ \ $h^{2,2}_{BC}$ \ \ &\ \ \ $h^{0,1}_{\bar{\partial}}+h^{2,0}_{\bar{\partial}}  -k^{2,0}$ \\ 
& \  &  \  &  \   & \   \\ 
$p=3$ & $k^{3,0}$  &  $h^{0,2}_{\bar{\partial}}  $  &  $h^{0,1}_{\bar{\partial}}+h^{2,0}_{\bar{\partial}} -k^{2,0} $   & 1    
\end{tabular}
\end{center}

where  $h^{3,1}_2 = dim_{\mathbb{C}} (E^{3,1}_2)$, 
and $h^{2,2}_{BC}$ is given by:
\begin{eqnarray*}
h^{2,2}_{BC} &=& - {\check{h}}^1({\mathcal{PH}})  +h^{0,1}_{\bar{\partial}} -h^{0,2}_{\bar{\partial}} -h^{1,0}_{\bar{\partial}}+ h^{1,1}_{\bar{\partial}} \\
\  & \ & \ \ \ \ \ \ \ + h^{1,2}_{\bar{\partial}}  +h^{2,0}_{\bar{\partial}}   +2 h^{3,1}_2 +k^{1,1}-k^{1,2}-2k^{2,0}  \  .
\end{eqnarray*}
We also note that by a result of Tosatti{\cite{Tosatti} which is covered later in the article we have
$$h^{1,1}_{BC}={\check{h}}^1({\mathcal{PH}}) = 2h^{0,1}_{\bar{\partial}} - b^1 +dim_{\mathbf{R}}(H^{1,1}({\mathbf{R}}))  $$
where $H^{1,1}({\mathbf{R}})$ is the subgroup of $H^2({\mathbf{R}})$ of deRham classes which have a representative which is a real $d$-closed $1,1$-form.  
\ \\

In section \ref{SeqOfCohMaps}, we derive some useful ``almost exact" sequences of cohomology on compact complex mandifolds of any dimension that prove useful in the calculation of Bott-Chern/Aeppli cohomology.  In section \ref{GeneralResults}, we derive Bott-Chern/Aeppli cohomology hodge numbers, $h^{1,0}_{BC}$, $h^{n-1,0}_{BC}$, and $h^{n,0}_{BC}$ on compact complex manifolds.   In section \ref{MoreGeneralResults}, we derive Bott-Chern/Aeppli cohomology hodge numbers, $h^{p,0}_{BC}$ and $h^{p,n}_{BC}$ on compact complex manifolds.   In section \ref{DecompositionOfH_A}, we give a decomposition of $H^{p,q}_{A}$ into subspaces $H^{p,q}_{\partial}/im({\bar{\partial}})$ $G^{p,q}_{\partial}$, and $L^{p,q}_{\partial}$ on compact complex manifolds that proves crucial in our calculation of $H^{p,q}_{BC}$ on compact complex 3-folds.  In section \ref{CompactComplex3folds}, we first derive a formula for $h^{1,1}_{BC}$ on a compact complex manifold.  We then restrict to compact complex 3-folds, to calculate the complete Bott-Chern/Aeppli cohomology, $H^{p,q}_{BC}$.  In section \ref{Examples} we show how our table of Bott-Chern/Aeppli cohomology on a compact complex 3-fold applies to a hypothetical complex structure on $S^6$ and on how it is in agreement with Angella and Tomassini's calculation on a Calabi-Eckman 3-fold\cite{AngellaAndTomassini} and with Angella's calculation for small Kuranishi deformations of the Iwasawa manifold\cite{Angella}.
\ \\

I would like to thank Daniele Angella for making very helpful suggestions, insightful corrections and generous encouragement. 
I would also like to thank the organizers and hosts of the International Conference in Current Developments and New Directions in Kahler Geometry held at the University of Notre Dame in June, 2017.  In addition, I am grateful for support of this research by funds from the School of Science, Engineering and Technology, Pennsylvania State University-Harrisburg.  I thank Kevin Li and Ron Walker for help in revising the manuscript for publication.

\section{Sequences of maps of cohomology} \label{SeqOfCohMaps}
We first note that there is a natural projection map, $\pi_{\bar{\partial}} :H^{p,q}_{BC} \rightarrow H^{p,q}_{\bar{\partial}} $ where for a $d$-closed $p,q$-form, $\theta$, we have
$ \theta / im( \partial {\bar{\partial}})  \mapsto  \theta / im( {\bar{\partial}})  \  . $  Similiarly, we have the projection map,  $\pi_A :H^{p,q}_{\bar{\partial}} \rightarrow H^{p,q}_A $ where for a ${\bar{\partial}}$-closed $p,q$-form, $\theta$, we have
$ \theta / im(  {\bar{\partial}})  \mapsto  \theta /(im(\partial) + im( {\bar{\partial}}))  \  . $ 

Consider the following sequences of maps of cohomology on a compact $n$-dimensional complex manifold $X$ with, $p=0,\cdots,n$,
$$ 0 \rightarrow  H^{p,0}_{BC} \ \ \stackrel{\pi_{\bar{\partial}}}{\longrightarrow}  H^{p,0}_{\bar{\partial}} \ \ \stackrel{  {\tiny{\pi_A }}}{ \longrightarrow}  H^{p,0}_A  
\ \  \stackrel{{\bar{\partial}}}{\rightarrow} H^{p,1}_{BC} \ \  \stackrel{\pi_{\bar{\partial}}}{\longrightarrow}  \cdots $$
$$ \cdots \rightarrow  H^{p,q}_{BC} \ \ \stackrel{\pi_{\bar{\partial}}}{\longrightarrow}  H^{p,q}_{\bar{\partial}} \ \ \stackrel{  {\tiny{\pi_A }}}{ \longrightarrow}  H^{p,q}_A  
\ \  \stackrel{{\bar{\partial}}}{\rightarrow} H^{p,q+1}_{BC} \ \  \stackrel{\pi_{\bar{\partial}}}{\longrightarrow}  \cdots $$
$$ \stackrel{  {\tiny{\pi_A }}}{ \longrightarrow}  H^{p,n-1}_A  \ \  \stackrel{{\bar{\partial}}}{\rightarrow}  H^{p,n}_{BC} \ \  \stackrel{\pi_{\bar{\partial}}}{\longrightarrow}   H^{p,n}_{\bar{\partial}} \ \  \stackrel{  {\tiny{\pi_A }}}{ \longrightarrow}  H^{p,n}_A   \rightarrow  0  \ \  . $$
   
The above sequence of maps for the case of $X$, a hypothetical complex 3-fold diffeomorphic to $S^6$, was also given in McHugh{\cite{McHugh}}. It is pretty straight forward to show 
  the sequence of maps above is exact at $H^{p,q}_{BC}$.  Namely, 
$$ ker(\pi_{\bar{\partial}}  : H^{p,q}_{BC} \rightarrow  H^{p,q}_{\bar{\partial}} )  =  im( {\bar{\partial}} : H^{p,q-1}_A  \rightarrow H^{p,q}_{BC} )  \ \  .  $$
It is also straight forward to show the sequence of maps above is exact at $H^{p,q}_A$.  Namely, 
$$im( \pi_A  : H^{p,q}_{\bar{\partial}} \rightarrow  H^{p,q}_A  )  = ker( {\bar{\partial}} : H^{p,q}_A  \rightarrow H^{p,q+1}_{BC} )  \ \ \ .  $$
We have  the following ``almost" exactness at $ H^{p,q}_{\bar{\partial}}$:
\begin{proposition}   
$$dim(Im( \pi_{\bar{\partial}}: H^{p,q}_{BC} \rightarrow H^{p,q}_{\bar{\partial}} ) )   
 =   k^{p,q}  +  dim(ker(  \pi_A :  H^{p,q}_{\bar{\partial}} \rightarrow  H^{p,q}_A ) ) $$
\end{proposition}
{\bf{Proof}}:  We first show that 
$$ker(  \pi_A :  H^{p,q}_{\bar{\partial}} \rightarrow  H^{p,q}_A )  \subseteq Im( \pi_{\bar{\partial}}: H^{p,q}_{BC} \rightarrow H^{p,q}_{\bar{\partial}} ) \  . $$ 
Indeed, if we have a ${\bar{\partial}}$-closed $p,q$-form, $\phi$ such that $\phi = \partial \lambda + {\bar{\partial}} \chi$, then $\partial \lambda$ is ${\bar{\partial}}$-closed and represents the same element as $\phi$ in $H^{p,q}_{\bar{\partial}} $.  Clearly, $\partial \lambda \in ker(\partial) \cap ker({\bar{\partial}})$ and thus
 $$[\phi ] = [ \partial \lambda ] \in Im( \pi_{\bar{\partial}}: H^{p,q}_{BC} \rightarrow H^{p,q}_{\bar{\partial}} ) \ . $$  Denote $ V^{p,q} = Im( \pi_{\bar{\partial}}: H^{p,q}_{BC} \rightarrow H^{p,q}_{\bar{\partial}} )$.  
We also have by what we have just shown that 
$$ker(  \pi_A :  H^{p,q}_{\bar{\partial}} \rightarrow  H^{p,q}_A ) = ker(  \pi_A :  V^{p,q} \rightarrow  H^{p,q}_A ) \ . $$
If we apply the Rank Theorem from basic Linear Algebra, to $ \pi_A :  V^{p,q} \rightarrow  H^{p,q}_A $ we get the result of our proposition: 
$$dim(V^{p,q})  = 
 dim(ker( \pi_A :  V^{p,q} \rightarrow  H^{p,q}_A )) 
 + dim(im(\pi_A :  V^{p,q} \rightarrow  H^{p,q}_A )) $$
or \\
$dim(Im( \pi_{\bar{\partial}}: H^{p,q}_{BC} \rightarrow H^{p,q}_{\bar{\partial}} ) )   = $
$$ \ \ \ \ \ \ \ \ \ \ \ \ \ \ \ \ \ \ \  \ \ \ \ \ \ \ \ \ \ \ \ \ \ \ \  \  dim(ker(  \pi_A :  H^{p,q}_{\bar{\partial}} \rightarrow  H^{p,q}_A ) ) 
+ dim(\frac{ker( \partial  ) \cap ker( {\bar{\partial}} ) }{im( \partial  ) + im( {\bar{\partial}}) }) \  . $$
{\bf{End of Proof}} 

Using the Rank theorem from basic Linear Algebra we have for $0 \leq p \leq n$ and $ 0 \leq q \leq n $  :  
 \begin{eqnarray*}
 h^{p,q}_{BC}  &  =  &  dim(Ker( \pi_{\bar{\partial}}: H^{p,q}_{BC} \rightarrow H^{p,q}_{\bar{\partial}} ))+ dim(Im( \pi_{\bar{\partial}}: H^{p,q}_{BC} \rightarrow H^{p,q}_{\bar{\partial}} ))  \ , \\
 k^{p,q} + h^{p,q}_{\bar{\partial}}   &  =  &  k^{p,q} + dim(Ker(  \pi_A :  H^{p,q}_{\bar{\partial}} \rightarrow  H^{p,q}_A ) )    \\ 
\ & \ & \ \ \ \ \ \ \ \ \ \ \ \ \ \ \ \ \ \  \ \ \ \ \ \ \ \ \ \ \ \ \ \ \ \  + \ dim(Im(\pi_A  : H^{p,q}_{\bar{\partial}} \rightarrow  H^{p,q}_A  ) )   \ , \\ 
h^{p,q}_A    &  =  &  dim( ker( {\bar{\partial}} : H^{p,q}_A  \rightarrow H^{p,q+1}_{BC} ))  + dim(im( {\bar{\partial}} : H^{p,q}_A  \rightarrow H^{p,q+1}_{BC} ) )  \  . \\
\end{eqnarray*}
Notice that we have trivially added $k^{p,q}$ to both sides of the equation from the Rank theorem for $h^{p,q}_{\bar{\partial}}$.
Creating an alternating sum over $q$ of both sides of the equations above and using the lemmas and proposition on exactness, we have:
\begin{eqnarray}\label{AlmostExactSequenceTheorem}
 \sum^{n}_{q=0} (-1)^q ( h^{p,q}_{BC} - (h^{p,q}_{\bar{\partial}} + k^{p,q}) + h^{p,q}_A)  = 0
\end{eqnarray} or using the duality, $h^{p,q}_A = h^{n-p,n-q}_{BC}$,
$$ \sum^{n}_{q=0} (-1)^q ( h^{p,q}_{BC} + h^{n-p,n-q}_{BC}) = \sum^{n}_{q=0} (-1)^{q} (h^{p,q}_{\bar{\partial}} + k^{p,q}) \ \  .$$
Using the Riemann-Roch-Hirzebruch Theorem, we can also write this equation as
\begin{eqnarray*}
 \sum^{n}_{q=0} (-1)^q ( h^{p,q}_{BC} + h^{n-p,n-q}_{BC}-k^{p,q}) & =  & \sum^{n}_{q=0} (-1)^{q} h^{p,q}_{\bar{\partial}} \\
  \  & = & \chi(\Omega^p)  \\
  \  & = &  \int_X ch(\Omega^p)td(X) 
  \end{eqnarray*}
where $\Omega^p$ is the bundle of holomorphic $p,0$-forms.  Thus, 
\begin{proposition}\label{rrhbc}
On an $n$-dimensional compact complex manifold, $X$,
$$ \sum^{n}_{q=0} (-1)^q ( h^{p,q}_{BC} + h^{n-p,n-q}_{BC}-k^{p,q}) =  \int_X ch(\Omega^p)td(X) $$
is a ``topological" invariant of our compact complex manifold (in that it depends only on the topological structure of $\Omega^p$ and $TX$).
\end{proposition}
We shall use the (almost exact) sequences above in calculating the Bott-Chern cohomology.  In more specific situations the sequences split into two or more useful sequences.

\section{Some general results for Aeppli/Bott-Chern Cohomology} \label{GeneralResults}
We have the following straight forward result (see McHugh\cite{McHugh}) when $b^1=0$
\begin{proposition}
If $b^1 = 0$ then $H^{1,0}_{BC} = H^{0,1}_{BC}=H^{n-1,n}_A=H^{n,n-1}_A = 0$
\end{proposition}
{\bf Proof}  Let $[\phi] \in H^{1,0}_{BC}$ where $\phi$ is a $1,0$-form such that $\partial \phi = 0$ and ${\bar{\partial}}=0$.  Specifically, $d \phi = 0$ and so since $b^1=0$ we have 
$\phi = d f$ for some global function $f$.  Now $\phi = \partial f + {\bar{\partial}} f $.   Since $\phi$ is a $1,0$-form, we have ${\bar{\partial}} f =0$.  Thus, $f$ is a global holomorphic function on a compact complex manifold and thus must be a constant function.  Hence,  $\partial f = 0$, and we have $\phi =0$. \\

We can generalize this to the case $b^1 \neq 0$.  Consider the following portion of our sequence of maps:
$$0 \rightarrow H^{0,0}_{BC} \rightarrow H^{0,0}_{\bar{\partial}} \rightarrow  H^{0,0}_{A} \rightarrow H^{0,1}_{BC} \rightarrow  H^{0,1}_{\bar{\partial}} \rightarrow\ldots   \ \  \  \  \  . $$
For a (connected) compact complex manifold, $b^0 = h_{\bar{\partial}}^{0,0}=k^{0,0}=1$, and it is easy to see that $h_{BC}^{0,0} =1$.  Thus the image of the map, $H^{0,0}_{\bar{\partial}} \rightarrow H^{0,0}_A$ is one-dimensional.  Also, $H^{0,0}_A$ consists of all pluri-harmonic scalar functions on our compact manifold.  These must be constant by the maximum modulus theorem so $h^{0,0}_A=1$. ( I thank Michael Albanese for pointing this out to me.\cite{Albanese}) By exactness, the kernel of the map, $H^{0,1}_{BC} \rightarrow H^{0,1}_{\bar{\partial}}$ is zero. Since, 
\begin{eqnarray*}
K^{0,1} &=& (ker(d)\cap {\mathcal{E}}^{0,1})/(Im(\partial )\cap {\mathcal{E}}^{0,1} + Im({\bar{\partial}})\cap {\mathcal{E}}^{0,1})  \\
\ & = &   (ker(d)\cap {\mathcal{E}}^{0,1})/ (Im({\bar{\partial}})\cap {\mathcal{E}}^{0,1}) \ \ \  ,
\end{eqnarray*}
 we have, 
\begin{proposition} On a compact complex manifold with $b^1 \neq 0$, $h^{0,1}_{BC} = k^{0,1} $ .
\end{proposition}  
We will also show the following:
\begin{proposition}
For any $n$-dimensional complex manifold, $X$,  $h^{n-1,0}_{BC}=h^{n-1,0}_{\bar{\partial}} \  . $
\end{proposition} 
Before we prove this we need the following lemma:
 \begin{lemma} For any $n$-dimensional complex manifold, $X$, a ${\bar{\partial}}$-closed, ${\partial}$-exact $n,0$-form is zero.  Explicitly, let $\phi$ be an $n-1,0$-form.
 If $\partial \phi$ is ${\bar{\partial}}$-closed, then $\partial \phi = 0$.
 \end{lemma}
 {\bf Proof of Lemma} Following in the same manner almost verbatim as Lemma 2.2 in Brown\cite{Brown},
\begin{eqnarray*}
 \int_X \partial \phi \wedge {\bar{\partial}}{\bar{ \phi}} & =& \int_X {\bar{\partial}}(\partial \phi \wedge {\bar{\phi}})  - \int_X ({\bar{\partial}}(\partial \phi ))\wedge {\bar{\phi}} \\
 \  &=& \int_X d(\partial \phi \wedge {\bar{\phi}})  - \int_X ({\bar{\partial}}(\partial \phi ))\wedge {\bar{\phi}} = 0+0 =0 
\end{eqnarray*}
by Stokes theorem.  

Locally, $\partial \phi = f dz^1 \wedge \cdots \wedge dz^n$
and so 
\begin{eqnarray*}
\int_X \partial \phi \wedge {\bar{\partial}}{\bar{ \phi}} & = & \int_X |f|^2 dz^1 \wedge \cdots \wedge dz^n \wedge d{\bar{z}}^1 \wedge \cdots \wedge d{\bar{z}}^n \\
\ & =& \int_X (-1)^{\frac{n(n-1)}{2}} |f|^2 dx^1 \wedge dy^1 \wedge \cdots \wedge dx^n \wedge dy^n = 0
\end{eqnarray*}
Thus $f = 0$ and $ \partial \phi = 0$. \\
 {\bf Proof of Proposition}
We have the sequence,
$$0 \rightarrow H^{n-1,0}_{BC} \rightarrow H^{n-1,0}_{\bar{\partial}} \rightarrow H^{n-1,0}_A \rightarrow \cdots $$
which is exact at $H^{n-1,0}_{BC}$.   Thus $H^{n-1,0}_{BC} \hookrightarrow H^{n-1,0}_{\bar{\partial}} $ injectively.
We shall show that the map is surjective.  

Let $[\phi] \in H^{n-1,0}_{\bar{\partial}} $, with of course, ${\bar{\partial}} \phi =0 $.  We have that $\partial \phi $ is a ${\bar{\partial}}$-closed, ${\partial}$-exact $n,0$-form
and thus $\partial \phi = 0$.  
We have $ \phi \in ker(\partial) \cap ker({\bar{\partial}})$ and $[\phi ]_{\bar{\partial}} = \pi_{\bar{\partial}}( [\phi]_{BC} ) $.
Hence $H^{n-1,0}_{BC} = H^{n-1,0}_{\bar{\partial}} $ and $h^{n-1,0}_{BC}=h^{n-1,0}_{\bar{\partial}}$. \\

We shall also show in almost exactly the same way that $$h^{n,0}_{BC} = h^{n,0}_{\bar{\partial}} = k^{n,0} \ \  . $$  We repeat the argument for completeness. See also Angella\cite{Angella}, Section 1.4.4  .
\begin{proposition} \label{hBCn0}
For any $n$-dimensional complex manifold, $X$,
$$h^{n,0}_{BC}=h^{n,0}_{\bar{\partial}} = k^{n,0} \  . $$
\end{proposition}
{\bf Proof of Proposition}
We have the sequence,
$$0 \rightarrow H^{n,0}_{BC} \rightarrow H^{n,0}_{\bar{\partial}} \rightarrow H^{n,0}_A \rightarrow \cdots $$
which is exact at $H^{n,0}_{BC}$.   Thus $H^{n,0}_{BC} \hookrightarrow H^{n,0}_{\bar{\partial}} $ injectively.
We shall show that the map is surjective.  

Let $[\phi] \in H^{n,0}_{\bar{\partial}} $, with of course, ${\bar{\partial}} \phi =0 $.  We have that $\partial \phi =0 $ since $\phi$ is an $n,0$-form.
We have $ \phi \in ker(\partial) \cap ker({\bar{\partial}})$ and $[\phi ]_{\bar{\partial}} =\pi_{\bar{\partial}}( [\phi]_{BC} )$.
Hence $H^{n,0}_{BC} = H^{n,0}_{\bar{\partial}} $ and $h^{n,0}_{BC}=h^{n,0}_{\bar{\partial}}$. \\

 Define the projection map, $\pi_K :H^{p,q}_{BC} \rightarrow K $ where for a $d$-closed $p,q$-form, $\theta$, we have
$ \theta / im(  {\bar{\partial}})  \mapsto  \theta /(im(\partial) + im( {\bar{\partial}}))  $.  Now consider a nonzero, $[\phi]_{BC} \in H^{n,0}_{BC}$.  By our lemma above, $\phi$ is not $\partial$-exact.  Thus the map, 
$ \pi_K  :  H^{n,0}_{BC} \rightarrow K^{n,0} $
is injective and $h^{n,0}_{BC} \leq k^{n,0}$.  Combining this together we get
$h^{n,0}_{BC} \leq k^{n,0} \leq h^{n,0}_{\bar{\partial}} = h^{n,0}_{BC}$
and thus $h^{n,0}_{BC} = h^{n,0}_{\bar{\partial}} = k^{n,0}$.

\section{Some more general results for Aeppli/Bott-Chern Cohomology: Calculating $h^{p,0}_{BC}$ and $h^{p,0}_A$} \label{MoreGeneralResults}
Recall the following results just shown above. \\
On a compact complex manifold with $b^1 = 0$, $H^{1,0}_{BC} = H^{0,1}_{BC}=H^{n-1,n}_A=H^{n,n-1}_A = 0 $  .  \\
On a  compact complex manifold with $b^1 \neq 0$, $h^{0,1}_{BC} = k^{0,1}  $. \\
For any $n$-dimensional complex manifold, $X$, $h^{n-1,0}_{BC}=h^{n-1,0}_{\bar{\partial}}  $ and $h^{n,0}_{BC}=h^{n,0}_{\bar{\partial}} = k^{n,0} $. \\

We will be deriving the ``extreme" rows/columns of our Bott-Chern and Aeppli cohomology tables, i.e. $h^{p,0}_{BC}$ and $h^{p,0}_A$ for all $p$.
We first derive a formula for $h^{p,0}_{BC}$.   We have the beginning of the  sequence,
$$ 0 \rightarrow H^{p,0}_{BC} \rightarrow H^{p,0}_{\bar{\partial}} \rightarrow H^{p,0}_A  \rightarrow \cdots \ \ \ .$$
 We can consider the sequence,  
$$ 0 \rightarrow H^{p,0}_{BC} \rightarrow H^{p,0}_{\bar{\partial}}  \rightarrow Im( \pi_A : H^{p,0}_{\bar{\partial}}\rightarrow  H^{p,0}_A) \rightarrow 0  $$
which is exact except in the middle with the non-exactness measured by $K^{p,0}$.  Thus
$$ h^{p,0}_{BC} - (h^{p,0}_{\bar{\partial}} + k^{p,o}) +  dim( Im(\pi_A)): H^{p,0}_{\bar{\partial}}\rightarrow  H^{p,0}_A)) = 0 $$
and 
\begin{eqnarray} \label{hp0_BC}
  h^{p,0}_{BC} =h^{p,0}_{\bar{\partial}}  -  (dim( Im( \pi_A: H^{p,0}_{\bar{\partial}}\rightarrow  H^{p,0}_A)) -k^{p,0} ) \ .  
\end{eqnarray}
We claim that 
$$ dim( Im(\pi_A : H^{p,0}_{\bar{\partial}}\rightarrow  H^{p,0}_A)) - k^{p,0}   = dim( Im(\partial :H^{p,0}_{\bar{\partial}} \rightarrow H^{p+1,0}_{\bar{\partial}})) $$
where $\partial :H^{p,0}_{\bar{\partial}} \rightarrow H^{p+1,0}_{\bar{\partial}}$ is the map from the Frolicher spectral sequence.
Note that we may consider the embedding, 
$$K^{p,0} \hookrightarrow  Im(\pi_A: H^{p,0}_{\bar{\partial}}\rightarrow  H^{p,0}_A)) $$ and thus 
$$dim( Im( \pi_A: H^{p,0}_{\bar{\partial}}\rightarrow  H^{p,0}_A)) - k^{p,0}  = dim( Im( \pi_A: H^{p,0}_{\bar{\partial}}\rightarrow  H^{p,0}_A)/ K^{p,0}) \ \ . $$
We will thus show that $$dim( Im( \pi_A: H^{p,0}_{\bar{\partial}}\rightarrow  H^{p,0}_A)/ K^{p,0}) = dim( Im(\partial :H^{p,0}_{\bar{\partial}} \rightarrow H^{p+1,0}_{\bar{\partial}})) \ \ .$$
 Consider an  non-zero element,
 $$ [\theta]_A/K^{p,0} \in   Im( \pi_A: H^{p,0}_{\bar{\partial}}\rightarrow  H^{p,0}_A)/ K^{p,0}$$  
 for  $\theta$, a nonzero ${\bar{\partial}}$-closed, $p,0$-form.  We will show that $\partial [\theta]_{\bar{\partial}} \neq 0$.  We can not have $\partial \theta = {\bar{\partial}} \nu $ as this is a $p+1,0$ form and if we have $\partial \theta = 0$ then $[\theta]_A \in K^{p,0}$ contradicting $ [\theta]_A/K^{p,0} \neq 0$.  Thus 
$$ \partial [\theta]_{\bar{\partial}} \in   Im(\partial :H^{p,0}_{\bar{\partial}} \rightarrow H^{p+1,0}_{\bar{\partial}})$$
 and is nonzero.  This is thus an injective linear map from  
$$ Im( \pi_A: H^{p,0}_{\bar{\partial}}\rightarrow  H^{p,0}_A)/ K^{p,0} \ {\mathrm{to}} \  Im(\partial :H^{p,0}_{\bar{\partial}} \rightarrow H^{p+1,0}_{\bar{\partial}}) \ . $$  We still need to show it is onto.   Clearly, the map is well defined.

Now consider a nonzero element,
$$\partial [\phi]_{\bar{\partial}} \in  Im(\partial :H^{p,0}_{\bar{\partial}} \rightarrow H^{p+1,0}_{\bar{\partial}}) $$
 for  $\phi$, a nonzero ${\bar{\partial}}$-closed, $p,0$-form.  If $ [\phi]_A /K^{p,0} = 0$, then $\phi = \kappa + \partial \mu$ where
$\kappa$ is a local representative of $K^{p,0}$.  But then $\partial \phi = 0$ and this is a contradiction.  Hence, $[\phi]_A /K^{p,0}$ is non-zero and thus $\partial [\phi]_{\bar{\partial}}$ is in the image of our map just above.  Our map is onto.    Thus we have an isomorphism,
$$ Im(\pi_A: H^{p,0}_{\bar{\partial}}\rightarrow  H^{p,0}_A)/ K^{p,0}   = Im(\partial :H^{p,0}_{\bar{\partial}} \rightarrow H^{p+1,0}_{\bar{\partial}}) $$
and
$$ {\mathrm{dim}}( Im( \pi_A: H^{p,0}_{\bar{\partial}}\rightarrow  H^{p,0}_A)) - k^{p,0}   \  \  = {\mathrm{dim}}( Im(\partial :H^{p,0}_{\bar{\partial}} \rightarrow H^{p+1,0}_{\bar{\partial}})) \ . $$

We now calculate $ {\mathrm{dim}}( Im(\partial :H^{p,0}_{\bar{\partial}} \rightarrow H^{p+1,0}_{\bar{\partial}})) $.
Recall the Frolicher Sequence,
$$0 \rightarrow H^{0,0}_{\bar{\partial}} \stackrel{\partial}{\rightarrow} \cdots  \stackrel{\partial}{\rightarrow} H^{p,0}_{\bar{\partial}} \stackrel{\partial}{\rightarrow} H^{p+1,0}_{\bar{\partial}} \stackrel{\partial}{\rightarrow} \cdots   \stackrel{\partial}{\rightarrow} H^{n,0}_{\bar{\partial}}  \stackrel{\partial}{\rightarrow} 0 $$
where the non-exactness is measured by the spectral sequence terms,
$$E^{q,0}_2 = Ker( H^{q,0}_{\bar{\partial}} \stackrel{\partial}{\rightarrow} H^{q+1,0}_{\bar{\partial}})/ Im( H^{q-1,0}_{\bar{\partial}} \stackrel{\partial}{\rightarrow} H^{q,0}_{\bar{\partial}} )\ . $$
Let $ker_q = dim( Ker( H^{q,0}_{\bar{\partial}} \stackrel{\partial}{\rightarrow} H^{q+1,0}_{\bar{\partial}}))$ and $im_q = dim( Im( H^{q-1,0}_{\bar{\partial}} \stackrel{\partial}{\rightarrow} H^{q,0}_{\bar{\partial}}))$ and thus $h^{q,0}_2 = ker_q - im_q $.
Using, 
$$h^{q,0}_{\bar{\partial}} = ker_q +im_{q+1} = im_q +h_2^{q,0} +im_{q+1}$$
we have,
$$ h^{q,0}_{\bar{\partial}}- h_2^{q,0} =  im_q  +im_{q+1} \ . $$
Thus, performing a telescoping alternating sum,
$$im_{p+1}=( h^{p,0}_{\bar{\partial}} - h^{p,0}_2) - (h^{p-1,0}_{\bar{\partial}}-h^{p-1,0}_2) 
				+\ \ \ \ \ \ \ \ \ \ \ \ \ \ \ $$ $$ \ \ \ \ \ \ \ \ \ \ \ \ \ \ \ \ \ \ \ \ \ \ \ \ \ \ \  \\ \ \ \  \cdots +(-1)^{p-1}( h^{1,0}_{\bar{\partial}}-h_2^{1,0}) + (-1)^p( h^{0,0}_{\bar{\partial}}-h_2^{0,0})  \  .  $$
For compact complex manifolds, $h^{0,0}_{\bar{\partial}}=h^{0,0}_2 = 1 $ , so we can write
\begin{eqnarray}  \label{im_pPlus1}
 im_{p+1} = \sum_{j=1}^{p}(-1)^{p+j} ( h^{j,0}_{\bar{\partial}}-h_2^{j,0}) \ . 
\end{eqnarray}
Finally, using equation \ref{hp0_BC}, we have
$$  h^{p,0}_{BC} =h^{p,0}_{\bar{\partial}}  -  \sum_{j=1}^{p}(-1)^{p+j} ( h^{j,0}_{\bar{\partial}}-h_2^{j,0}) \ . $$

We now derive formulas for $h^{p,0}_A=h^{n,n-p}_{BC}$.  We will proceed by induction on $q$ for calculating $h^{n, q}_{BC}$.  We start the induction by deriving $h^{n, 1}_{BC}$.  We know from Proposition \ref{hBCn0} that $h^{n,0}_A = h^{0,n}_{BC} =  h^{n,0}_{BC}=k^{n,0}$.  
Thus the map, ${\bar{\partial}} : H^{n,0}_A \rightarrow H^{n,1}_{BC}$ is zero and we have the clipped sequence,
$$ 0 \rightarrow H^{n,1}_{BC} \rightarrow H^{n,1}_{\bar{\partial}} \rightarrow H^{n,1}_A \rightarrow H^{n,2}_{BC} \rightarrow \cdots $$
Hence, $H^{n,1}_{BC}$ embeds into $ H^{n,1}_{\bar{\partial}}$.
We have 
$$dim( H^{n,1}_{\bar{\partial}}) =dim( Ker(\pi_A: H^{n,1}_{\bar{\partial}}\rightarrow  H^{n,1}_A)) + dim(Im( \pi_A: H^{n,1}_{\bar{\partial}}\rightarrow  H^{n,1}_A))$$
so that
$$h^{n,1}_{\bar{\partial}} = dim(Im(\pi_{\bar{\partial}})) -k^{n,1} + dim(  Im(\pi_A: H^{n,1}_{\bar{\partial}}\rightarrow  H^{n,1}_A))$$
$$= h^{n,1}_{BC} + dim( Im(\pi_A: H^{n,1}_{\bar{\partial}}\rightarrow  H^{n,1}_A)) - k^{n,1} \ \ . $$
We claim that 
$$ dim(  Im( \pi_A: H^{n,1}_{\bar{\partial}}\rightarrow  H^{n,1}_A)) - k^{n,1} = 0$$ 
Note that we may consider $$K^{n,1} \hookrightarrow   Im(\pi_A: H^{n,1}_{\bar{\partial}}\rightarrow  H^{n,1}_A)) $$ and thus 
$$dim( Im(\pi_A: H^{n,1}_{\bar{\partial}}\rightarrow  H^{n,1}_A)) - k^{n,1}  = dim(  Im(\pi_A: H^{n,1}_{\bar{\partial}}\rightarrow  H^{n,1}_A)/ K^{n,1}) \ \ . $$
Indeed, consider an element,
 $$ [\theta]_A/K^{n,1} \in    Im(\pi_A: H^{n,1}_{\bar{\partial}}\rightarrow  H^{n,1}_A)/ K^{n,1}$$  
 for  $\theta$, a nonzero ${\bar{\partial}}$-closed, but not ${\bar{\partial}}$-exact, $n,1$-form. We automatically have that $\partial \theta =0$ as it is an $n,1$-form.
Thus $d\theta =0$ and  $[\theta]_A \in K^{n,1}$ with $ \theta = \kappa + \partial \mu + {\bar{\partial}} \nu $
and $ [\theta]_A/K^{n,1} =0$. \\  Hence 
$$ dim( Im( \pi_A: H^{n,1}_{\bar{\partial}}\rightarrow  H^{n,1}_A)/ K^{n,1}) =0$$
and
$$ {\mathrm{dim}}( Im(\pi_A: H^{n,1}_{\bar{\partial}}\rightarrow  H^{n,1}_A)) - k^{n,1}  = 0 \  \ . $$
Thus 
\begin{eqnarray} \label{hbc_n1}
 h^{0,n-1}_A=h^{n,1}_{BC} = h^{n,1}_{\bar{\partial}}   \  .  
\end{eqnarray} 
We shall now proceed with our induction on $q$ in $h^{n,q}_{BC}$.  Assume that we can write formulas for $h^{n,j}_{BC}$ in terms of $h^{*,*}_{\bar{\partial}}$,  $h^{*,*}_2$, and $k^{*,*}$, for 
$0 \leq j \leq q$.  We note that we already can write formulas for $h^{n,j}_A$, $0 \leq j \leq n $, in such terms since $h^{n,j}_A = h^{0,n-j}_{BC} = h^{n-j,0}_{BC}$.
Consider again the sequence,
$$ 0 \rightarrow H^{n,1}_{BC} \rightarrow \cdots \rightarrow H^{n,q}_{BC} \rightarrow H^{n,q}_{\bar{\partial}} \rightarrow H^{n,q}_A \ \ \ \ \ \ \ \ \ \ \ \ \ \ \ \ \ \ \ \ $$ $$ \ \ \ \ \ \ \ \ \ \ \ \ \ \ \ \ \ \ \ \ \ \ \ \ \ \ \ \ \ \ \ \ \ \ \ \  \rightarrow H^{n,q+1}_{BC} \rightarrow  H^{n,q+1}_{\bar{\partial}} \rightarrow H^{n,q+1}_A \rightarrow  \cdots  $$
and shorten it to 
$$0 \rightarrow H^{n,1}_{BC} \rightarrow \cdots \rightarrow H^{n,q}_{BC} \rightarrow H^{n,q}_{\bar{\partial}} \rightarrow H^{n,q}_A  \ \ \ \ \ \ \ \ \ \ \ \ \ \ \ \ \ \ \ \ $$ $$ \ \ \ \ \ \ \ \ \ \ \ \ \ \ \ \ \ \ \ \ \ \ \ \ \ \ \  \rightarrow  H^{n,q+1}_{BC} \rightarrow  H^{n,q+1}_{\bar{\partial}} \rightarrow   Im( \pi_A: H^{n,q+1}_{\bar{\partial}}\rightarrow  H^{n,q+1}_A) \rightarrow  0 \ . $$
We shall show that $  Im(\pi_A: H^{n,q+1}_{\bar{\partial}}\rightarrow  H^{n,q+1}_A)) = K^{n,q+1}$.  It is clear that 
$$ K^{n,q+1} \subseteq Im(\pi_A: H^{n,q+1}_{\bar{\partial}}\rightarrow  H^{n,q+1}_A) \ . $$
If $[\phi ]_A \in  Im(\pi_A: H^{n,q+1}_{\bar{\partial}}\rightarrow  H^{n,q+1}_A)$
 then ${\bar{\partial}} \phi =0 $ and $\partial \phi = 0$ since $\phi$ is an $n,q+1$-form.
Thus $[\phi ]_A \in K^{n,q+1}$ and  
 $$ Im(\pi_A: H^{n,q+1}_{\bar{\partial}}\rightarrow  H^{n,q+1}_A) \subseteq K^{n,q+1} \ . $$
Thus
$$ Im(\pi_A: H^{n,q+1}_{\bar{\partial}}\rightarrow  H^{n,q+1}_A) = K^{n,q+1} \ . $$
We have then that 

$$ {\sum_{l=1}^{q} (-1)^{l-1}(h^{n,l}_{BC} - (h^{n,l}_{\bar{\partial}} + k^{n,l}) + h^{n,l}_A ) }\ \ \ \ \ \ \ \ \ \ \ \ \ \ \ \   \ \ \ \ \ \ \ \ \ \ \ \ \ \ \ \ \ \ \ \ $$ $$ \ \ \ \ \ \ \ \ \ \ \ \ \ \ \ \ \ \ \ \ \ \ \ \ \ \ \ \ \ \ \ \ \ \ \ \ \ \ \ \ \ \ \ \ \ \ \ \ + \  (-1)^{q}(h^{n,q+1}_{BC} - (h^{n,q+1}_{\bar{\partial}} + k^{n,q+1}) +k^{n,q+1})   = 0$$
and hence, 
\begin{eqnarray}\label{hbc_nqplus1}
h^{n,q+1}_{BC} & = & h^{n,q+1}_{\bar{\partial}} +\sum_{l=1}^{q} (-1)^{q+l}(h^{n,l}_{BC} - (h^{n,l}_{\bar{\partial}} + k^{n,l}) + h^{n,l}_A ) \end{eqnarray} 
 and 
\begin{eqnarray}\label{h_A_nqplus1}
h^{n-p,0}_A & = &  h^{0,p}_{\bar{\partial}} +\sum_{l=1}^{p-1} (-1)^{p-1+l}(h^{n-l,0}_{A} - (h^{0,n-l}_{\bar{\partial}} + k^{n,l}) + h^{n-l,0}_{BC} ) \  . 
\end{eqnarray} 
We now know $h^{p,0}_{BC}$ and $h^{p,0}_A$ in terms of $h^{*,*}_{\bar{\partial}}$,  $h^{*,*}_2$, and $k^{*,*}$ for all $p$ such that $0 \leq p \leq n$. 

\section{Decomposition of $H^{p,q}_A$} \label{DecompositionOfH_A}
In this section we define some vector spaces associated with $H^{p,q}_A$ and prove some results with regard to these
spaces.  These will be helpful in some later sections in calculating further Bott-Chern/Aeppli cohomology of compact complex 3-folds. 

 Consider the map: $ \partial : H^{p,q}_A \rightarrow H^{p+1,q}_{\bar{\partial}}  $ where $ [\theta ]_A \mapsto [\partial \theta ]_{\bar{\partial}}$ .
We define
$${\mathcal{G}}^{p,q}_{\partial} = ker(\partial:  H^{p,q}_A \rightarrow H^{p+1,q}_{\bar{\partial}}) \ , $$
$$ G^{p,q}_{\partial} = {\mathcal{G}}^{p,q}_{\partial} / Im( \ \pi_A: H^{p,q}_{\partial} \rightarrow H^{p,q}_A ) \ , $$
and
$$ L^{p,q}_{\partial} = H^{p,q}_A/{\mathcal{G}}^{p,q}_{\partial} \ . $$
Here, $\pi_A$ is the obvious projection map from $ H^{p,q}_{\partial}$ to $H^{p,q}_A$  exactly similar to $ \pi_A: H^{p,q}_{\bar{\partial}} \rightarrow H^{p,q}_A $.
Following terminology of Popovici{\cite{Popovici}}, we suggest calling $ G^{p,q}_{\partial}$, {\em strongly Gauduchon cohomology} and calling $ L^{p,q}_{\partial}$, {\em weakly Gauduchon cohomology}.  (Popovici{\cite{Popovici}} calls an hermitian metric on a complex n-fold, {\em strongly Gauduchon}, if the $n-1$ power of its associated $1,1$-form, $\omega$, is such that $\partial (\omega^{n-1}) = {\bar{\partial}} \eta$ for some $n,n-2$-form, $\eta$.)
We also define the hodge numbers,  $g^{p,q}_{\partial} = dim_{\mathbf{C}} (G^{a,b}_{\partial})$ and $l^{p,q}_{\partial} = dim_{\mathbf{C}} (L^{a,b}_{\partial})$. \\  \  \\
We define  ${\mathcal{G}}^{p,q}_{\bar{\partial}}$, $G^{p,q}_{\bar{\partial}}$,  $L^{p,q}_{\bar{\partial}}$, $g^{p,q}_{\bar{\partial}}$, and $l^{p,q}_{\bar{\partial}}$ completely analagously:
$${\mathcal{G}}^{p,q}_{\bar{\partial}} = ker({\bar{\partial}}:  H^{p,q}_A \rightarrow H^{p,q+1}_{\partial}) \ , $$
$$ G^{p,q}_{\bar{\partial}} = {\mathcal{G}}^{p,q}_{\bar{\partial}} / Im( \ \pi_A : H^{p,q}_{\bar{\partial}} \rightarrow H^{p,q}_A ) \ , $$
and
$$ L^{p,q}_{\bar{\partial}} = H^{p,q}_A/{\mathcal{G}}^{p,q}_{\bar{\partial}} \  $$
with also $g^{p,q}_{\bar{\partial}} = dim_{\mathbf{C}} (G^{a,b}_{\bar{\partial}})$ and $l^{p,q}_{\bar{\partial}} = dim_{\mathbf{C}} (L^{a,b}_{\bar{\partial}})$.
We note that $g^{p,q}_{\partial} = g^{q,p}_{\bar{\partial}}$, $l^{p,q}_{\partial} = l^{q,p}_{\bar{\partial}}$, $g^{p,q}_{\partial} = g^{p+1,q-1}_{\bar{\partial}}$
and $  l^{n-p-2,n-q}_{\partial} = h^{p+1,q}_{\bar{\partial}} - l^{p,q}_\partial - k^{p+1,q}$.  The first two equations are obvious from the definitions.  We shall give a proof of the last equation later on.  We give a proof of the third equation now.
\begin{proposition}  \label{gpartialiandgbarpartial} $$g^{p,q}_{\partial} = g^{p+1,q-1}_{\bar{\partial}}$$
\end{proposition}
We shall show an isomorphism, $\phi$, between $G^{p,q}_{\partial}$ and  $G^{p+1,q-1}_{\bar{\partial}}$.  Indeed, let $[\mu]_{  G^{p,q}_{\partial}} \in G^{p,q}_{\partial}$, where $\mu$ is a $p,q$-form such that $\partial \mu = {\bar{\partial}} \nu$ for some $p+1,q-1$-form $\nu$.  Note that we may consider $[\nu]_{  G^{p+1,q-1}_{\bar{\partial}}} \in G^{p+1,q-1}_{\bar{\partial}}$.  Thus we will show
$$\phi : G^{p,q}_{\partial} \rightarrow G^{p+1,q-1}_{\bar{\partial}} $$
 $$\phi(\, [\mu]_{  G^{p,q}_{\partial}}\, ) = [\nu]_{  G^{p+1,q-1}_{\bar{\partial}}}  $$
is a well defined and bijective linear map.  
Let us proceed to show that it is well defined.

If $\tilde{\mu}$ is another $p,q$-form such that $[{\tilde{\mu}}]_{  G^{p,q}_{\partial}} = [\mu]_{  G^{p,q}_{\partial}}$, then 
$\tilde{\mu} =  \mu + \chi + \partial \sigma + {\bar{\partial}}\tau $ 
where $\partial \chi = 0$.  Thus,
$\partial \tilde{\mu} = {\bar{\partial}}( \nu  - \partial \tau) $ 
and 
$$\phi( \, [{\tilde{\mu}}]_{  G^{p,q}_{\partial}}\, ) = [\nu -\partial \tau ]_{  G^{p+1,q-1}_{\bar{\partial}}} =  [\nu]_{  G^{p+1,q-1}_{\bar{\partial}}}  $$
Our map, $\phi$ is clearly linear.
Now to show $\phi$ is one-to-one.  If $ [{\hat{\mu}}]_{  G^{p,q}_{\partial}} $ is such that 
$\phi( \, [{\hat{\mu}}]_{  G^{p,q}_{\partial}} \, ) =  [\nu]_{  G^{p+1,q-1}_{\bar{\partial}}}  $
then $$ \partial \mu -\partial {\hat{\mu}} = {\bar{\partial}} \nu  -  {\bar{\partial}} ( \nu  + \partial \tau ) \  \  .  $$
Thus  $ {\hat{\mu}}  = \mu  + \chi  - {\bar{\partial}} \tau $
for some $p,q$-form, $\chi$ such that $\partial \chi =0$
and
$ [{\hat{\mu}}]_{  G^{p,q}_{\partial}} =  [ \mu ]_{  G^{p,q}_{\partial}}  \  .  $
Thus our map is one-to-one.  

To show our map is onto, let  
$[\nu]_{  G^{p+1,q-1}_{\bar{\partial}}} \in  G^{p+1,q-1}_{\bar{\partial}}  $.  Then
$ {\bar{\partial}} \nu = \partial \mu $
for some $p,q$-form $\mu$.  We see that  
$[\mu]_{  G^{p,q}_{\partial}} \in G^{p,q}_{\partial}$
and that 
$ \phi( \, [\mu]_{  G^{p,q}_{\partial}} \, )  = [\nu]_{  G^{p+1,q-1}_{\bar{\partial}}} $. 
Thus our map is onto.  Thus we have that  $G^{p,q}_{\partial}$ is isomorphic to  $G^{p+1,q-1}_{\bar{\partial}}$.
\ \\

Underlying much of our analysis will be the following decompositions of $H^{p,q}_{BC}$ and $H^{p,q}_A$ (dependent on some, always existing, choice of hermitian metric):
\begin{lemma}\label{decomp}
$$H^{p,q}_{BC} = \partial G^{p-1,q} \oplus \partial L^{p-1,q} \oplus {\bar{\partial}} L^{p,q-1} \oplus K^{p,q} \ \ , $$
$$H^{p,q}_{BC} ={\bar{\partial}} G^{p,q-1} \oplus {\bar{\partial}} L^{p,q-1} \oplus \partial L^{p-1,q} \oplus K^{p,q} \ \ , $$
$$H^{p,q}_A = H^{p,q}_\partial / ({\bar{\partial}} L^{p,q-1}_{\bar{\partial}} ) \oplus G^{p,q}_\partial \oplus L^{p,q}_\partial \ \ , $$
$$H^{p,q}_A = H^{p,q}_{\bar{\partial}} / (\partial L^{p,q-1}_\partial ) \oplus G^{p,q}_{\bar{\partial}} \oplus L^{p,q}_{\bar{\partial}} \ \ , $$
\end{lemma}
  The proof of the second two statements is straightforward from definitions. We give a proof of the first two statements at least with respect to the hodge numbers:
\ \\ 
We notice that 
$$dim_{\mathbf{C}} ( im( {\bar{\partial}} : H^{p,q}_A \rightarrow H^{p,q+1}_{BC})) = g^{p,q}_{\bar{\partial}} + l^{p,q}_{\bar{\partial}} $$
since 
$$ h^{p,q}_A = dim_{\mathbf{C}} ( ker( {\bar{\partial}} : H^{p,q}_A \rightarrow H^{p,q+1}_{BC}))  + dim_{\mathbf{C}} ( im( {\bar{\partial}} : H^{p,q}_A \rightarrow H^{p,q+1}_{BC})) \  , $$
$$dim_{\mathbf{C}} ( ker( {\bar{\partial}} : H^{p,q}_A \rightarrow H^{p,q+1}_{BC})) = dim_{\mathbf{C}} ( im( H^{p,q}_{\bar{\partial}} \rightarrow H^{p,q}_{A}))  $$
and 
$$ h^{p,q}_A =  dim_{\mathbf{C}} ( im( H^{p,q}_{\bar{\partial}} \rightarrow H^{p,q}_{A})) + g^{p,q}_{\bar{\partial}} + l^{p,q}_{\bar{\partial}} \  . $$
We also know that $$ l^{p-1,q+1}_\partial +k^{p,q+1} = dim_{\mathbf{C}} ( im( \pi_{\bar{\partial}}: H^{p,q+1}_{BC} \rightarrow H^{p,q+1}_{\bar{\partial}}))  \  . $$

By the exactness of our sequence at $H^{p,q+1}_{BC}$ we have
$$dim_{\mathbf{C}} ( ker(\  \pi_{\bar{\partial}}: H^{p,q+1}_{BC} \rightarrow H^{p,q+1}_{\bar{\partial}})  = g^{p,q}_{\bar{\partial}} + l^{p,q}_{\bar{\partial}} $$
and by the rank theorem,
\begin{eqnarray*}
h^{p,q+1}_{BC} &=& g^{p,q}_{\bar{\partial}} + l^{p,q}_{\bar{\partial}} +  l^{p-1,q+1}_\partial +k^{p,q+1} \\
\ & = & g^{p-1,q+1}_\partial + l^{p-1,q+1}_\partial + l^{p,q}_{\bar{\partial}} +k^{p,q+1} \  .
\end{eqnarray*}
or
$$h^{p,q}_{BC} = g^{p,q-1}_{\bar{\partial}} + l^{p,q-1}_{\bar{\partial}} +  l^{p-1,q}_\partial +k^{p,q}$$ and
$$h^{p,q}_{BC} = g^{p-1,q}_\partial + l^{p-1,q}_\partial + l^{p,q-1}_{\bar{\partial}} +k^{p,q} \ .$$
We also give here the formulas from the decompositions for $H^{p,q}_A$:
$$h^{p,q}_A = h^{p,q}_\partial - l^{p,q-1}_{\bar{\partial}} + g^{p,q}_\partial +l^{p,q}_\partial$$
and
$$h^{p,q}_A = h^{p,q}_{\bar{\partial}} - l^{p,q-1}_\partial+ g^{p,q}_{\bar{\partial}} +l^{p,q}_{\bar{\partial}} \ .$$
 \  \\
We have the following result which gives the formula: 
\begin{lemma}     For $0 \leq p \leq n-2$ and $0 \leq q \leq n$, we have: \\
   $$l^{n-p-2,n-q}_\partial =  h^{p+1,q}_{\bar{\partial}} -l^{p,q}_\partial - k^{p+1,q} \  \  .$$
\end{lemma}  
{\bf Proof} Consider the portion of our sequence 
$$ \ldots \rightarrow H^{p+1,q-1}_A \rightarrow H^{p+1,q}_{BC} \rightarrow H^{p+1,q}_{\bar{\partial}} \rightarrow H^{p+1,q}_A \rightarrow \ldots $$
We have by our exactness results
$$Im(H^{p+1,q}_{BC} \rightarrow H^{p+1,q}_{\bar{\partial}}) = \partial L^{p,q}_\partial \oplus K^{p+1,q} \ . $$
Thus $$0 \rightarrow H^{p+1,0}_{BC} \rightarrow \ldots \rightarrow H^{p+1,q-1}_A \rightarrow H^{p+1,q}_{BC} \rightarrow  \partial L^{p,q}_\partial \oplus K^{p+1,q} \rightarrow 0   $$
and $$l^{p,q}_{\partial}+k^{p+1,q} = h^{p+1,q}_{BC}-h^{p+1,q-1}_A + \ldots + (-1)^{\epsilon (q)} h^{p+1,0}_{BC} \ . $$
In a similiar manner we consider the portion of our sequence:
$$ \ldots \rightarrow H^{n-p-1,n-q}_{BC} \rightarrow H^{n-p-1,n-q}_{\bar{\partial}} \rightarrow H^{n-p-1,n-q}_A \rightarrow \ldots $$
and with $Ker( H^{n-p-1,n-q}_{\bar{\partial}} \rightarrow H^{n-p-1,n-q}_A) = \partial L^{n-p-2, n-q}_\partial $ we have
$$ 0 \rightarrow \partial L^{n-p-2,n-q}_{\partial} \rightarrow H^{n-p-1,n-q}_{\bar{\partial}} \rightarrow H^{n-p-1,n-q}_A \rightarrow \ldots \rightarrow H^{n-p-1,n}_A \rightarrow 0 \ . $$
Thus,
\begin{eqnarray*}
 l^{n-p-2,n-q}_{\partial} & = & h^{n-p-1,n-q}_{\bar{\partial}} - h^{n-p-1,n-q}_A + \ldots + (-1)^{\epsilon (q)+1} h^{n-p-1,n}_A \\
\   & = & h^{p+1,q}_{\bar{\partial}} - h^{p+1,q}_{BC} + \ldots + (-1)^{\epsilon (q)+1} h^{p+1,0}_{BC} 
\end{eqnarray*}
Thus
$$  l^{n-p-2,n-q}_{\partial} = h^{p+1,q}_{\bar{\partial}} - l^{p,q}_\partial - k^{p+1,q} \ . $$
\  \\

We need a further decomposition of $H^{p,q}_A$ that we develop here.   Denote $H^{p,q}_{\bar{\partial}} / im(\partial) = {\mathcal{H}}^{p,q}_{\bar{\partial}}$ and similiarly,  $H^{p,q}_{\partial} / im({\bar{\partial}}) = {\mathcal{H}}^{p,q}_{\partial}$.  
Define $$ {\mathcal{H}}^{p,q}_{\bar{\partial}} \sqcap  {\mathcal{H}}^{p,q}_{\partial} =  {\mathcal{H}}^{p,q}_{\bar{\partial}} \cap  {\mathcal{H}}^{p,q}_\partial = K^{p,q} \ \  .  $$
We also have the short filtration:
$$ {\mathcal{H}}^{p,q}_{\bar{\partial}} \sqcap  {\mathcal{H}}^{p,q}_\partial \subseteq   {\mathcal{H}}^{p,q}_{\bar{\partial}} \cap {\mathcal{G}}^{p,q}_\partial \subseteq  {\mathcal{H}}^{p,q}_{\bar{\partial}} \ \  .$$  
We also define $$  {\mathcal{H}}^{p,q}_{\bar{\partial}} \sqcap G^{p,q}_\partial  = ( {\mathcal{H}}^{p,q}_{\bar{\partial}} \cap {\mathcal{G}}^{p,q}_\partial ) / ( {\mathcal{H}}^{p,q}_{\bar{\partial}} \sqcap  {\mathcal{H}}^{p,q}_\partial ) $$
 and 
$$  {\mathcal{H}}^{p,q}_{\bar{\partial}} \sqcap L^{p,q}_\partial  = (   {\mathcal{H}}^{p,q}_{\bar{\partial}} ) / (  {\mathcal{H}}^{p,q}_{\bar{\partial}} \cap {\mathcal{G}}^{p,q}_\partial ) \  . $$
Thus, we have the decomposition,
$$  {\mathcal{H}}^{p,q}_{\bar{\partial}} = ( {\mathcal{H}}^{p,q}_{\bar{\partial}} \sqcap  {\mathcal{H}}^{p,q}_\partial ) \  \oplus \  ({\mathcal{H}}^{p,q}_{\bar{\partial}} \sqcap G^{p,q}_\partial ) \  \oplus \   ({\mathcal{H}}^{p,q}_{\bar{\partial}}\sqcap L^{p,q}_\partial ) \  . $$
\ \\

Consider the filtration,
$$ {\mathcal{H}}^{p,q}_{\bar{\partial}} \cap  {\mathcal{H}}^{p,q}_\partial \subseteq {\mathcal{G}}^{p,q}_{\bar{\partial}} \cap  {\mathcal{H}}^{p,q}_\partial \subseteq  {\mathcal{G}}^{p,q}_{\bar{\partial}} \cap  {\mathcal{G}}^{p,q}_\partial \subseteq {\mathcal{G}}^{p,q}_{\bar{\partial}} \  . $$
We can define 
$$ {\mathcal{G}}^{p,q}_{\bar{\partial}} \sqcap  {\mathcal{H}}^{p,q}_\partial  =  ( {\mathcal{G}}^{p,q}_{\bar{\partial}} \cap  {\mathcal{H}}^{p,q}_\partial ) \  ,  $$
$$ G^{p,q}_{\bar{\partial}} \sqcap  {\mathcal{H}}^{p,q}_\partial  =  ( {\mathcal{G}}^{p,q}_{\bar{\partial}} \cap  {\mathcal{H}}^{p,q}_\partial ) / ( {\mathcal{H}}^{p,q}_{\bar{\partial}} \cap  {\mathcal{H}}^{p,q}_\partial ) $$
and 
$$ {\mathcal{G}}^{p,q}_{\bar{\partial}} \sqcap  G^{p,q}_\partial  =  ( {\mathcal{G}}^{p,q}_{\bar{\partial}} \cap  {\mathcal{G}}^{p,q}_\partial ) / ( {\mathcal{G}}^{p,q}_{\bar{\partial}} \cap  {\mathcal{H}}^{p,q}_\partial ) \  .  $$
Recall that 
$$ {\mathcal{H}}^{p,q}_{\bar{\partial}} \sqcap  G^{p,q}_\partial  =  ( {\mathcal{H}}^{p,q}_{\bar{\partial}} \cap  {\mathcal{G}}^{p,q}_\partial ) / ( {\mathcal{H}}^{p,q}_{\bar{\partial}} \cap  {\mathcal{H}}^{p,q}_\partial )  $$
and so we define
$$ G^{p,q}_{\bar{\partial}} \sqcap  G^{p,q}_\partial  =  ( {\mathcal{G}}^{p,q}_{\bar{\partial}} \sqcap  G^{p,q}_\partial ) / ( {\mathcal{H}}^{p,q}_{\bar{\partial}} \sqcap  G^{p,q}_\partial ) \  .  $$
We also define
$$ {\mathcal{G}}^{p,q}_{\bar{\partial}} \sqcap  L^{p,q}_\partial =  {\mathcal{G}}^{p,q}_{\bar{\partial}} /  ( {\mathcal{G}}^{p,q}_{\bar{\partial}} \cap  {\mathcal{G}}^{p,q}_\partial ) \  .  $$
Recall that 
$$ {\mathcal{H}}^{p,q}_{\bar{\partial}} \sqcap  L^{p,q}_\partial  =   {\mathcal{H}}^{p,q}_{\bar{\partial}}  / ( {\mathcal{H}}^{p,q}_{\bar{\partial}} \cap  {\mathcal{G}}^{p,q}_\partial )  $$
and so we define
$$ G^{p,q}_{\bar{\partial}} \sqcap  L^{p,q}_\partial  =  ( {\mathcal{G}}^{p,q}_{\bar{\partial}} \sqcap  L^{p,q}_\partial ) / ( {\mathcal{H}}^{p,q}_{\bar{\partial}} \sqcap  L^{p,q}_\partial ) \  .  $$
Thus we can write the decomposition,
$$  G^{p,q}_{\bar{\partial}} = ( G^{p,q}_{\bar{\partial}} \sqcap  {\mathcal{H}}^{p,q}_\partial ) \  \oplus \  ( G^{p,q}_{\bar{\partial}} \sqcap G^{p,q}_\partial ) \  \oplus \   (G^{p,q}_{\bar{\partial}}\sqcap L^{p,q}_\partial ) \  . $$
\ \\

Continuing, we define
$$  L^{p,q}_{\bar{\partial}} \sqcap  {\mathcal{H}}^{p,q}_\partial  = (   {\mathcal{H}}^{p,q}_\partial ) / (  {\mathcal{G}}^{p,q}_{\bar{\partial}} \cap {\mathcal{H}}^{p,q}_\partial ) \  ,  $$
$$ L^{p,q}_{\bar{\partial}} \sqcap  {\mathcal{G}}^{p,q}_\partial =  {\mathcal{G}}^{p,q}_\partial /  ( {\mathcal{G}}^{p,q}_{\bar{\partial}} \cap  {\mathcal{G}}^{p,q}_\partial ) \  ,  $$
$$ L^{p,q}_{\bar{\partial}} \sqcap  G^{p,q}_\partial =  ( L^{p,q}_{\bar{\partial}} \sqcap  {\mathcal{G}}^{p,q}_\partial ) /  ( L^{p,q}_{\bar{\partial}} \sqcap  {\mathcal{H}}^{p,q}_\partial ) \  ,  $$
and 
$$ L^{p,q}_{\bar{\partial}} \sqcap  L^{p,q}_\partial =  ( L^{p,q}_{\bar{\partial}}  ) /  ( L^{p,q}_{\bar{\partial}} \sqcap  {\mathcal{G}}^{p,q}_\partial ) \  .   $$
Thus we have the decomposition
$$  L^{p,q}_{\bar{\partial}} = ( L^{p,q}_{\bar{\partial}} \sqcap  {\mathcal{H}}^{p,q}_\partial ) \  \oplus \  ( L^{p,q}_{\bar{\partial}} \sqcap G^{p,q}_\partial ) \  \oplus \   (L^{p,q}_{\bar{\partial}}\sqcap L^{p,q}_\partial ) \  . $$
\  \\

We  will use the following:
\begin{lemma} \label{ker}
$$Ker(\partial: H^{p,q}_{\bar{\partial}} \rightarrow H^{p+1,q}_{\bar{\partial}}) = G^{p,q}_\partial  \sqcap \, ( H^{p,q}_{\bar{\partial}} / im(\partial)) \  \oplus \  K^{p,q} \  \oplus \  \partial L^{p-1,q}_\partial  \ . $$
\end{lemma}
{\bf{Proof}}  Recall from Lemma \ref{decomp} the decompsoition of $H^{p,q}_{BC}$ :
$$H^{p,q}_{BC} ={\bar{\partial}} G^{p,q-1} \oplus {\bar{\partial}} L^{p,q-1} \oplus \partial L^{p-1,q} \oplus K^{p,q} \ \ , $$
and thus that 
$$ Im( \, \pi_{\bar{\partial }} \, : \, H^{p,q}_{BC} \rightarrow H^{p,q}_{\bar{\partial}} ) \cong K^{p,q} \oplus  \partial L^{p-1,q} \  .  $$
We see that 
$$ Im( \, \pi_{\bar{\partial }} \, : \, H^{p,q}_{BC} \rightarrow H^{p,q}_{\bar{\partial}} ) \subseteq Ker(\partial: H^{p,q}_{\bar{\partial}} \rightarrow H^{p+1,q}_{\bar{\partial}}) 	\  .  $$
We also have that 
\begin{eqnarray*}
H^{p,q}_{\bar{\partial}} &  \cong &  H^{p,q}_{\bar{\partial}} / im(\partial) \oplus \partial L^{p-1,q} \\
\  & \cong &  H^{p,q}_\partial / (im({\bar{\partial}}) ) \sqcap H^{p,q}_{\bar{\partial}} / im(\partial) \oplus G^{p,q}_\partial \sqcap H^{p,q}_{\bar{\partial}} / im(\partial) \ \ \ \ \ \ \ \ \ \   \\
\  & \ &  \ \ \ \ \ \ \ \ \ \  \ \ \ \ \  \ \ \ \ \ \ \ \ \ \ \ \ \ \ \ \ \ \ \ \  \oplus L^{p,q}_\partial \sqcap H^{p,q}_{\bar{\partial}} / im(\partial) \,  \oplus \, \partial L^{p-1,q} \\
\  & \cong &  K^{p,q}  \oplus G^{p,q}_\partial \sqcap H^{p,q}_{\bar{\partial}} / im(\partial)\oplus L^{p,q}_\partial \sqcap H^{p,q}_{\bar{\partial}} / im(\partial) \,  \oplus \, \partial L^{p-1,q} 
\end{eqnarray*}
where we have used the fact that 
$$ K^{p,q} =( H^{p,q}_\partial / im({\bar{\partial}}) ) \cap (H^{p,q}_{\bar{\partial}} / im(\partial))  \  . $$
By the definitions of $G^{p,q}_\partial$ and $L^{p,q}_\partial$ we see that 
$$Ker(\partial: H^{p,q}_{\bar{\partial}} \rightarrow H^{p+1,q}_{\bar{\partial}}) = G^{p,q}_\partial  \sqcap \, ( H^{p,q}_{\bar{\partial}} / im(\partial)) \  \oplus \  K^{p,q} \  \oplus \  \partial L^{p-1,q}_\partial  \ . $$

\section{BC-A cohomology on generic compact complex 3-folds}\label{CompactComplex3folds}
We shall now complete  a table for a generic compact complex 3-fold.  
Using the formulas above for $h^{p,0}_{BC}$ and $h^{p,n}_{BC}$, one obtains: 
\begin{center}
$h^{0,0}_{BC} = 1$,    $h^{1,0}_{BC} = k^{1,0}$,  $h^{2,0}_{BC} = h^{2,0}_{\bar{\partial}}$,  $h^{3,3}_{BC} = 1 $, \\     $h^{3,0}_{BC} = h^{3,0}_{\bar{\partial}}= k^{3,0}$, \\ 
 $h^{3,1}_{BC} = h^{3,1}_{\bar{\partial}}=h^{0,2}_{\bar{\partial}}$,   \\  $h^{3,2}_{BC}= h^{0,1}_{\bar{\partial}} +h^{2,0}_{\bar{\partial}} - k^{3,1}$ .      
\end{center}

Using the Bigolin resolution\cite{Bigolin}\cite{Wells} of the sheaf of pluri-harmonic functions, ${\mathcal{PH}}$, we can also derive a formula for $h^{1,1}_{BC}$ for any compact complex manifold:
\begin{theorem} For any compact complex manifold,
$ H^{1,1}_{BC}=  {\check{H}}^{1}({\mathcal{PH}}) \  . $
Thus, $h^{1,1}_{BC} = h^{n-1,n-1}_A = |{\check{H}}^{1}({\mathcal{PH}})| \  .$
\end{theorem}
{\bf{Proof}}:  Let ${\mathcal{E}}^{p,q}$ denote the sheaf of $C^\infty$ $p,q$-forms.  The Bigolin resolution for ${\mathcal{PH}}$,
$$0 \rightarrow \  {\mathcal{PH}} \hookrightarrow {\mathcal{E}}^{ 0,0} \  {\stackrel{\partial {\bar{\partial}}}{\rightarrow}} \  {\mathcal{E}}^{ 1,1} {\stackrel{d}{\rightarrow}} \ {\mathcal{E}}^{ 2,1} \oplus {\mathcal{E}}^{ 1,2} \rightarrow {\mathrm{Coker}}(d) \rightarrow 0 $$
is acyclic (except at the last term,  ${\mathrm{Coker}}(d) $).  Thus 
\begin{eqnarray*}
 {\check{H}}^{1}({\mathcal{PH}}) & = & { \frac{ ker(d:\Gamma ({\mathcal{E}}^{ 1,1}) \rightarrow \Gamma ({\mathcal{E}}^{ 2,1}\oplus {\mathcal{E}}^{ 1,2} ) )}{ im(\partial {\bar{\partial}}:\Gamma ({\mathcal{E}}^{ 0,0}) \rightarrow \Gamma ({\mathcal{E}}^{ 1,1}) )}} \\
 \   & = &  \frac{ ker(d:C^{\infty \, 1,1} \rightarrow C^{\infty \, 2,1}\oplus C^{\infty \, 1,2} )}{im(\partial {\bar{\partial}}:C^{\infty \, 0,0} \rightarrow C^{\infty \, 1,1} ) } \\
\  & = & H^{1,1}_{BC}
\end{eqnarray*}
Following Bigolin\cite{Bigolin}, we note that 
$${\mathcal{PH}} = {\mathcal{O}} + {\bar{\mathcal{O}}} $$
and that we have the short exact sequence
$$0 \rightarrow {\mathbf{C}} \rightarrow {\mathcal{O}} \oplus {\bar{\mathcal{O}}} \rightarrow  {\mathcal{O}} + {\bar{\mathcal{O}}} \rightarrow 0 $$
or 
$$ 0 \rightarrow {\mathbf{C}} \rightarrow {\mathcal{O}} \oplus {\bar{\mathcal{O}}} \rightarrow {\mathcal{PH}} \rightarrow 0 $$
Hence we have the following portion of the subsequent long exact sequence,
$$0 \rightarrow {\check{H}}^{0}({\mathbf{C}}) \rightarrow  {\check{H}}^{0}({\mathcal{O}} \oplus {\bar{\mathcal{O}}}) \rightarrow  {\check{H}}^{0}({\mathcal{PH}})   \rightarrow \ \ \ \ \ \ \ \ \ \ \ \ \ \ \ \ \ \ \  $$ $$ \ \ \ \ \ \ \ \ \ \ \ \ \ \ \ \ \ \ \ \ \ \ \ \ {\check{H}}^{1}({\mathbf{C}}) \rightarrow  {\check{H}}^{1}({\mathcal{O}} \oplus {\bar{\mathcal{O}}}) \rightarrow  {\check{H}}^{1}({\mathcal{PH}}) \rightarrow  {\check{H}}^{2}({\mathbf{C}}) \rightarrow \cdots $$
which on a compact complex manifold becomes (using the fact that global holomorphic, global anti-holomorphic and global pluri-harmonic functions are constant) 
$$0 \rightarrow {\mathbf{C}} \rightarrow  {\mathbf{C}} \oplus {\mathbf{C}} \rightarrow  {\mathbf{C}}   \rightarrow  {\check{H}}^{1}({\mathbf{C}}) \rightarrow  \ \ \ \ \ \ \ \ \ \ \ \ \ \ \ \ \ \ \  $$ $$ \ \ \ \ \ \ \ \ \ \ \ \ \ \ \ \ \ \ \ \ \ \ \ \  {\check{H}}^{1}({\mathcal{O}}) \oplus  {\check{H}}^{1}({\bar{\mathcal{O}}}) \rightarrow  {\check{H}}^{1}({\mathcal{PH}}) \rightarrow  {\check{H}}^{2}({\mathbf{C}}) \rightarrow \cdots   \  .$$
Exactness at the first three terms of the sequence gives that the map ${\mathbf{C}} \rightarrow  {\check{H}}^{1}({\mathbf{C}})$ must be the zero map. Hence we have the clipped long exact sequence,
$$0 \rightarrow  {\check{H}}^{1}({\mathbf{C}}) \rightarrow  {\check{H}}^{1}({\mathcal{O}}) \oplus  {\check{H}}^{1}({\bar{\mathcal{O}}}) \rightarrow  {\check{H}}^{1}({\mathcal{PH}}) \rightarrow  {\check{H}}^{2}({\mathbf{C}}) \rightarrow \cdots   \  .$$
We rewrite this as,
$$0 \rightarrow  H^{1}({\mathbf{C}}) \rightarrow  H^{0,1}_{\bar{\partial}} \oplus  H^{1,0}_\partial \rightarrow H^{1,1}_{BC} \rightarrow  H^{2}({\mathbf{C}}) \rightarrow \cdots   \  .$$
If $b^2=0$ then $H^{2}({\mathbf{C}}) =0$ and we have the short exact sequence,
$$0 \rightarrow  H^{1}({\mathbf{C}}) \rightarrow  H^{0,1}_{\bar{\partial}} \oplus  H^{1,0}_\partial \rightarrow H^{1,1}_{BC} \rightarrow  0  \  .$$
Thus,
\begin{theorem} On a compact complex manifold with $b^2=0$, we have
$$ h^{1,1}_{BC} = h^{n-1,n-1}_A = 2h^{0,1}_{\bar{\partial}} - b^1 \ .$$
\end{theorem}
Our formula for $h^{1,1}_{BC}$ generalizes even further by results of Tosatti\cite{Tosatti} {\em which we follow here virtually verbatim }.  Tosatti gives the short exact sequence of sheaves,
$$ 0 \rightarrow {\mathbf{R}} \rightarrow {\mathcal{O}} \stackrel{Im}{ \rightarrow } {\mathcal{P}} \rightarrow 0 $$
where ${\mathcal{P}}$ is the sheaf of {\em real valued} pluriharmonic functions.  He thus also gives the resulting, long exact sequence in sheaf cohomoloby,
$$0 \rightarrow H^{0}({\mathbf{R}}) \rightarrow H^{0}({\mathcal{O}})  \rightarrow H^{0}({\mathcal{P }}) \ \ \ \ \ \ \ \ \ \ \ \ \ \ \ \ \ \ \ $$
$$ \ \ \ \ \ \ \ \ \ \ \ \ \ \ \ \ \ \ \rightarrow H^{1}({\mathbf{R}}) \rightarrow H^{1}({\mathcal{O}})  \rightarrow H^{1}({\mathcal{P }})   \rightarrow H^{2}({\mathbf{R}}) \rightarrow \cdots $$
The first three terms form a short exact sequence,
$$0 \rightarrow {\mathbf{R}} \rightarrow  {\mathbf{R}}^2 \rightarrow  {\mathbf{R}}  \rightarrow 0 $$
so that one has,
$$ 0 \rightarrow H^{1}({\mathbf{R}}) \rightarrow H^{1}({\mathcal{O}})  \rightarrow H^{1}({\mathcal{P }})   \rightarrow H^{2}({\mathbf{R}}) \stackrel{\pi^{0,2}}{\rightarrow} H^{2}({\mathcal{O}}) \rightarrow \cdots \ . $$
The map, $\pi^{0,2}$, is the projection to the $0,2$ part of $H^2({\mathbf{R}}) \cong H^2_{deRham}({\mathbf{R}})$ followed by the isomorphism, $H^2({\mathcal{O}}) \cong H^{0,2}_{\bar{\partial}}$. We have $$ker(\pi^{0,2}) = H^{1,1}({\mathbf{R}}) $$ 
where $H^{1,1}({\mathbf{R}})$ is the subgroup of $H^2({\mathbf{R}})$ of deRham classes which have a representative which is a real $d$-closed $1,1$-form. 
Thus one has the exact sequence of {\em real} vector spaces,
$$0 \rightarrow H^{1}({\mathbf{R}}) \rightarrow H^{1}({\mathcal{O}}) \rightarrow H^1({\mathcal{P}}) \rightarrow H^{1,1}({\mathbf{R}}) \rightarrow 0$$
Note that 
\begin{eqnarray*}
dim_{\mathbf{R}}(H^{1}({\mathbf{R}})) & = & dim_{\mathbf{C}}(H^{1}({\mathbf{C}}))  =  b^1 \\
dim_{\mathbf{R}}(H^{1}({\mathcal{O}})) & = & 2 dim_{\mathbf{C}}(H^{1}({\mathcal{O}})) =  2h^{0,1}_{\bar{\partial}} 
\end{eqnarray*}
and
\begin{eqnarray*}
dim_{\mathbf{R}}(H^{1}({\mathcal{P}})) & = & dim_{\mathbf{C}}(H^{1}({\mathcal{PH}}))  = h^{1,1}_{BC} 
\end{eqnarray*}
Thus
$$ b^1 - 2h^{0,1}_{\bar{\partial}} + h^{1,1}_{BC} - dim_{\mathbf{R}}(H^{1,1}({\mathbf{R}})) =0$$
and
\begin{theorem} On a compact complex manifold, we have
$$h^{1,1}_{BC} = 2h^{0,1}_{\bar{\partial}} - b^1 +dim_{\mathbf{R}}(H^{1,1}({\mathbf{R}})) \ . $$
\end{theorem}
If we can calculate $h^{1,2}_{BC}$ or $h^{2.2}_{BC}$ then we can know the BC-A cohomology completely in terms of Dolbeault and Frolicher terms.
Consider the following expression of $h^{2,1}_A$,
\begin{eqnarray*}
h^{2,1}_A & = & h^{2,1}_\partial - l^{2,0}_{\bar{\partial}} + g^{2,1}_\partial + l^{2,1}_\partial \\
 \  & = & h^{1,2}_{\bar{\partial}} - l^{0,2}_\partial + g^{2,1}_\partial + l^{2,1}_\partial
\end{eqnarray*}
We know from the almost exact sequence
$$0 \rightarrow L^{2,1}_\partial \rightarrow H^{3,1}_{\overline{\partial}} \rightarrow H^{3,1}_A \rightarrow H^{3,2}_{BC} \rightarrow H^{3,2}_{\overline{\partial}} \rightarrow H^{3,2}_A \rightarrow 0 $$ (which we note is exact at each term  except at   $ H^{3,2}_{\overline{\partial}}$) that $$  l^{2,1}_\partial = h^{3,1}_{\bar{\partial}} - h^{3,1}_A+h^{3,2}_{BC} -(h^{3,2}_{\bar{\partial}}+k^{3,2})+h^{3,2}_A \ .$$
We have from equation \ref{hbc_nqplus1} and equation \ref{hbc_n1} that
$$h^{3,2}_{BC} = h^{3,2}_{\bar{\partial}} + (-1)^2 (h^{3,1}_{BC} -(h^{3,1}_{\bar{\partial}} +k^{3,1})+ h^{3,1}_A) \ \ \ \ \ \ \ \ \ $$
$$ \ \ \ \ \ \ \ \ \ \ \  = h^{0,1}_{\bar{\partial}} + h^{3,1}_{\bar{\partial}}- h^{3,1}_{\bar{\partial}}-k^{2,0} +h^{2,0}_{\bar{\partial}}  \ \ \ $$
$$ \ \ \ \ \ \ \ \ \ \ \   \ \ \ \ \ \ \ \ \ \ \ \ \  = h^{0,1}_{\bar{\partial}} -k^{2,0} +h^{2,0}_{\bar{\partial}} $$
Plugging this into the expression for $l^{2,1}_\partial$ just above, we have
 $$l^{2,1}_\partial = h^{0,2}_{\bar{\partial}} - h^{2,0}_{\bar{\partial}} +h^{0,1}_{\bar{\partial}}-k^{2,0}+h^{2,0}_{\bar{\partial}} -(h^{0,1}_{\bar{\partial}}+k^{3,2})+k^{1,0} $$ $$ = h^{0,2}_{\bar{\partial}} -k^{2,0} \ . $$
Notice that this calculation follows through without using the assumptions, $b^1 =0$ or $b^2=0$.   We have instead used
the fact that $k^{3,2} = k^{1,0}$, from a Serre duality that can be proved using harmonic representations of $H_{A}$ and $H_{BC}$:
 $$K^{p,q} = H^{p,q}_{\bar{\partial}}/im(\partial ) \cap H^{p,q}_\partial / im({\bar{\partial}})\subseteq H^{p,q}_A \ , $$ and thus 
$$ \star K^{p,q} = K^{n-p, n-q} \subseteq H^{n-p,n-q}_{BC} \ . $$
\ \\
We will now show $g^{2,1}_\partial = 0$.  We know that $$g^{2,1}_\partial = g^{3,0}_{\bar{\partial}} = g^{0,3}_\partial \ . $$
Consider the following expression of $h^{0,3}_A$,
\begin{eqnarray*}
h^{0,3}_A & = & h^{0,3}_\partial - l^{0,2}_{\bar{\partial}} + g^{0,3}_\partial + l^{0,3}_\partial \\
 \  & = & h^{0,3}_{\bar{\partial}} - l^{2,0}_\partial + g^{0,3}_\partial + l^{0,3}_\partial
\end{eqnarray*}
Previously, we proved 
 For any $n$-dimensional complex manifold, $X$, a ${\bar{\partial}}$-closed, ${\partial}$-exact $n,0$-form is zero.  Explicitly, let $\phi$ be an $n-1,0$-form. If $\partial \phi$ is ${\bar{\partial}}$-closed, then $\partial \phi = 0$.
In other words, for any compact complex manifold, $l^{n-1,0}_\partial =0$.  Thus,we have above, $l^{2,0}_\partial = 0$.
We also showed previously that $h^{0,3}_A = h^{0,3}_{\bar{\partial}}= k^{0,3}$.
Thus $g^{2,1}_\partial = g^{0,3}_\partial = 0 $ and $l^{0,3}_\partial = 0$.

We proceed to calculate $l^{0,2}_\partial$.
Recall that 
 $$l^{n-p-2,n-q}_\partial =  h^{p+1,q}_{\bar{\partial}} -l^{p,q}_\partial - k^{p+1,q} \ .$$
Thus
$$l^{0,2}_\partial = h^{2,1}_{\bar{\partial}} - l^{1,1}_\partial - k^{2,1} \ . $$
To calculate $ l^{1,1}_\partial$ we consider the equation,
$$Ker(\partial: H^{p,q}_{\bar{\partial}} \rightarrow H^{p+1,q}_{\bar{\partial}}) = G^{p,q}_\partial  \sqcap \, ( H^{p,q}_{\bar{\partial}} / im(\partial)) \  \oplus \  K^{p,q} \  \oplus \  \partial L^{p-1,q}_\partial   $$
with $p,q \, = \, 2,1$ .
Thus 
$$|Ker(\partial: H^{2,1}_{\bar{\partial}} \rightarrow H^{3,1}_{\bar{\partial}})| = k^{2,1} +|G^{2,1}_\partial  \sqcap \, ( H^{2,1}_{\bar{\partial}} / im(\partial))| + l^{1,1}_\partial \ . $$
We have $g^{2,1}_\partial =0$ so 
$$|G^{2,1}_\partial  \sqcap \, ( H^{2,1}_{\bar{\partial}} / im(\partial))| =0 \ . $$
Thus 
$$  l^{1,1}_\partial =|Ker(\partial: H^{2,1}_{\bar{\partial}} \rightarrow H^{3,1}_{\bar{\partial}})| - k^{2,1} \ . $$
Define as temporary shorthand,
$$ker_{p,q} = |Ker(\partial: H^{p,q}_{\bar{\partial}} \rightarrow H^{p+1,q}_{\bar{\partial}})| $$ so that we have,
$$  l^{1,1}_\partial =ker_{2,1} - k^{2,1} \ . $$
Also define as temporary shorthand,
$$im_{p,q} = |Im(\partial: H^{p-1,q}_{\bar{\partial}} \rightarrow H^{p,q}_{\bar{\partial}})| \ .  $$
We shall calculate $ker_{2,1}$ from the Frolicher spectral sequence as we did above for $im_{p,0}$.
We have
$$h^{p,q}_2 = ker_{p,q} - im_{p,q}$$ and
$$h^{p,q}_{\bar{\partial}} = ker_{p,q} +im_{p+1,q} = ker_{p,q} +ker_{p+1,q} -h^{p+1,q}_2 \ . $$
Thus,
\begin{eqnarray*}
h^{p,q}_{\bar{\partial}} &= & ker_{p,q} +ker_{p+1,q} -h^{p+1,q}_2  \\
h^{p+1,q}_{\bar{\partial}} & =&  ker_{p+1,q} +ker_{p+2,q} -h^{p+2,q}_2 \\
\ & \vdots & \  \\
h^{n-1,q}_{\bar{\partial}} &= & ker_{n-1,q} +ker_{n,q} -h^{n,q}_2  \\
h^{n,q}_{\bar{\partial}} &= & ker_{n,q} 
\end{eqnarray*}
and
$$\sum^{n}_{j=p}(-1)^{j-p} h^{j,q}_{\bar{\partial}} = ker_{p,q} - \sum^{n-1}_{j=p} (-1)^{j-p} h^{j+1,q}_2 \ . $$
Thus,
$$ker_{p,q} = \sum^{n}_{j=p}(-1)^{j-p} h^{j,q}_{\bar{\partial}} + \sum^{n-1}_{j=p} (-1)^{j-p} h^{j+1,q}_2 \ . $$
For $n=3$ and $p,q \, = \, 2,1$ this is simply,
$$ker_{2,1} = h^{2,1}_{\bar{\partial}} - h^{3,1}_{\bar{\partial}} + h^{3,1}_2 \ . $$
Hence,
$$l^{1,1}_\partial = h^{2,1}_{\bar{\partial}} - h^{3,1}_{\bar{\partial}} + h^{3,1}_2 -k^{2,1}$$
and 
\begin{eqnarray*}
l^{0,2}_\partial &= &h^{2,1}_{\bar{\partial}} - l^{1,1}_\partial - k^{2,1} \\
\ & = & h^{2,1}_{\bar{\partial}} - ( h^{2,1}_{\bar{\partial}} - h^{3,1}_{\bar{\partial}} + h^{3,1}_2 -k^{2,1}) - k^{2,1} \\\
\ & = &  h^{3,1}_{\bar{\partial}} - h^{3,1}_2  \\
\end{eqnarray*}
Finally,
\begin{eqnarray*}
h^{2,1}_A & = & h^{1,2}_{\bar{\partial}} - l^{0,2}_\partial + g^{2,1}_\partial + l^{2,1}_\partial \\
 \  & = & h^{1,2}_{\bar{\partial}} - (h^{3,1}_{\bar{\partial}} - h^{3,1}_2) + 0 + ( h^{0,2}_{\bar{\partial}}-k^{2,0}  ) \\
\ & = & h^{1,2}_{\bar{\partial}} + h^{3,1}_2 - k^{2,0}
\end{eqnarray*}
Thus 
\begin{proposition}\label{hbc12}
On a general compact complex three-fold,
$$h^{1,2}_{BC} = h^{1,2}_{\bar{\partial}}   + h^{3,1}_2 - k^{2,0}$$
\end{proposition}
Equation  \ref{AlmostExactSequenceTheorem} with  $n=3$ and $p=2$ gives
\begin{eqnarray*}
h^{2,2}_{BC} &=&- h^{2,0}_{BC} + (h^{2,0}_{\bar{\partial}} +k^{2,0}) -h^{2,0}_A+ h^{2,1}_{BC} - (h^{2,1}_{\bar{\partial}} +k^{2,1}) \\
 \  & \ & \ \ \ \ \ \ \ +h^{2,1}_A  + (h^{2,2}_{\bar{\partial}} +k^{2,2}) -h^{2,2}_A+ h^{2,3}_{BC} - (h^{2,3}_{\bar{\partial}} +k^{2,3}) +h^{2,3}_A \\
 &=&- h^{2,0}_{\bar{\partial}} + h^{2,0}_{\bar{\partial}} +k^{2,0} -h^{0,2}_{\bar{\partial}}+ h^{1,2}_{\bar{\partial}}  + h^{3,1}_2- k^{2,0}- h^{2,1}_{\bar{\partial}} -k^{2,1} \\
 \  & \ & \ \ \ \ \    + h^{1,2}_{\bar{\partial}}  + h^{3,1}_2 - k^{2,0} + h^{1,1}_{\bar{\partial}} +k^{1,1}-2h^{0,1}_{\bar{\partial}} + b^1  \\
 \  & \ & \ \ \ \  -dim_{\mathbf{R}}(H^{1,1}({\mathbf{R}}))+h^{0,1}_{\bar{\partial}}+h^{2,0}_{\bar{\partial}}  -k^{2,0} - h^{1,0}_{\bar{\partial}} -k^{1,0} +k^{1,0} \\
\ &=&h^{2,0}_{\bar{\partial}}+ h^{1,2}_{\bar{\partial}} -h^{0,2}_{\bar{\partial}}  + h^{1,1}_{\bar{\partial}} -h^{1,0}_{\bar{\partial}} -h^{0,1}_{\bar{\partial}}+ b^1 \\
 \  & \ & \ \ \ \ \ \ \   -dim_{\mathbf{R}}(H^{1,1}({\mathbf{R}}))  +2 h^{3,1}_2 -2k^{2,0}+k^{1,1}-k^{1,2}  \  .
\end{eqnarray*}
\ \\
So,
\begin{eqnarray}\label{h22bc}
h^{2,2}_{BC} &=& -h^{0,1}_{\bar{\partial}} -h^{0,2}_{\bar{\partial}} -h^{1,0}_{\bar{\partial}}+ h^{1,1}_{\bar{\partial}}+ h^{1,2}_{\bar{\partial}}   +h^{2,0}_{\bar{\partial}}  +2 h^{3,1}_2 +k^{1,1}
\end{eqnarray}
\begin{eqnarray*}
\ \ \ \ \ \ \ \ \ \ \ \  \  & \ & \ \ \ \ \ \ \ \ \ \ \ \ -k^{1,2}-2k^{2,0} +b^1 -dim_{\mathbf{R}}(H^{1,1}({\mathbf{R}})) \  .
\end{eqnarray*}
\ \\
We also note for future use that we can just as easily write the formula for $h^{2,2}_{BC}$ in terms of $h^{1,1}_{BC} = {\check{h}}^1({\mathcal{PH}})$, Dolbeault, Frohlicher and $k^{i,j}$ :
\begin{eqnarray*}
h^{2,2}_{BC} &=& - h^{1,1}_{BC} +h^{0,1}_{\bar{\partial}} -h^{0,2}_{\bar{\partial}} -h^{1,0}_{\bar{\partial}}+ h^{1,1}_{\bar{\partial}} \\
\ & \  & \ \ \ \ \ \ \ \ \  \ \ \ \ \ \ \ \ \  \ \ \ \ \ \ \ \ \  + h^{1,2}_{\bar{\partial}}  +h^{2,0}_{\bar{\partial}}   +2 h^{3,1}_2 +k^{1,1}-k^{1,2}-2k^{2,0}  \  .
\end{eqnarray*}
We summarize with a table of our  results.
 \begin{center}
\begin{tabular}{c|c|c|c|c|}
\multicolumn{5}{c}{\bf{Bott-Chern cohomology on a Compact Complex 3-fold}}  \\
\ $h^{p,q}_{BC}$ & $q= 0$ & $ q=1$ & $q=2$ & $ q=3$ \\
\hline 
$p=0$ & \  \  \  \   1 \ \  \ \  &  $k^{1,0}$  & $h^{2,0}_{\bar{\partial}} $ & $k^{3,0}$ \\ 
& \  &  \  &  \   & \   \\ 
$p=1$ & $k^{1,0}$ & \  \ \  \ \ $h^{1,1}_{BC}$ \ \  & $ h^{1,2}_{\bar{\partial}}  + h^{3,1}_2 - k^{2,0} $ &  $h^{0,2}_{\bar{\partial}} $ \\ 
& \  &  \  &  \   & \   \\ 
$p=2$ & $h^{2,0}_{\bar{\partial}} $ & $ h^{1,2}_{\bar{\partial}}  + h^{3,1}_2 - k^{2,0} $ &  \ \ $h^{2,2}_{BC}$ \ \ &\ \ \ $h^{0,1}_{\bar{\partial}}+h^{2,0}_{\bar{\partial}}  -k^{2,0}$ \\ 
& \  &  \  &  \   & \   \\ 
$p=3$ & $k^{3,0}$  &  $h^{0,2}_{\bar{\partial}}  $  &  $h^{0,1}_{\bar{\partial}}+h^{2,0}_{\bar{\partial}} -k^{2,0} $   & 1  \\ 
& \  &   &  \   & \   
\end{tabular}
\end{center}
where 
$$h^{1,1}_{BC}=  2h^{0,1}_{\bar{\partial}} - b^1 +dim_{\mathbf{R}}(H^{1,1}({\mathbf{R}}))$$ 
and $h^{2,2}_{BC}$ is given in terms of $b^1$, $dim_{\mathbf{R}}(H^{1,1}({\mathbf{R}}))$, Dolbeault, Frolicher and $k^{p,q}$  in equation \ref{h22bc} .

\section{Examples of Bott-Chern/Aeppli Cohomology on  Compact Complex 3-folds}\label{Examples}
\subsection{BC-A cohomology on complex $S^6$}
For a hypothetical complex structure on $S^6$, it has been shown by Brown \cite{Brown}  that $h^{3,1}_2 = h^{0,2}_2 $, 
by Gray\cite{Gray} that $h^{3,0}_{\bar{\partial}} = 0$ and by Ugarte\cite{Ugarte} that
$h^{0,1}_{\bar{\partial}} =1+ h^{0,2}_{\bar{\partial}}$
and that $h^{1,1}_{\bar{\partial}}=1+h^{1,2}_{\bar{\partial}}+h^{1,0}_{\bar{\partial}}-h^{2,0}_{\bar{\partial}} \ . $ \\
(We also note of course that since $b^j=0$ for $0 < j < 6$, we have $k^{p,q} =0$ for $0 < p+q < 6$ .)

So for complex $S^6$ we can write from our results just above, $h^{1,2}_{BC} = h^{1,2}_{\bar{\partial}}  + h^{0,2}_2 \ . $
We can also write from our results just above, $ h^{1,1}_{BC} = 2h^{0,1}_{\bar{\partial}} = 2+2h^{0,2}_{\bar{\partial}} $
and
\begin{eqnarray*}
h^{2,2}_{BC} &=&  - h^{1,1}_{BC} +h^{0,1}_{\bar{\partial}} -h^{0,2}_{\bar{\partial}} -h^{1,0}_{\bar{\partial}}+ h^{1,1}_{\bar{\partial}}+ h^{1,2}_{\bar{\partial}}  +h^{2,0}_{\bar{\partial}}  +2 h^{3,1}_2   \\
\   &=&  -  2-2h^{0,2}_{\bar{\partial}}  +1+h^{0,2}_{\bar{\partial}} -h^{0,2}_{\bar{\partial}} -h^{1,0}_{\bar{\partial}}+ h^{1,1}_{\bar{\partial}}+ h^{1,2}_{\bar{\partial}}  +h^{2,0}_{\bar{\partial}}  +2 h^{3,1}_2   \\
\  &=&  -  1-2h^{0,2}_{\bar{\partial}}   -h^{1,0}_{\bar{\partial}}+ (1+h^{1,2}_{\bar{\partial}}+h^{1,0}_{\bar{\partial}}-h^{2,0}_{\bar{\partial}})+ h^{1,2}_{\bar{\partial}}  +h^{2,0}_{\bar{\partial}}  +2 h^{0,2}_2   \\
\ & = & 2 h^{1,2}_{\bar{\partial}}  +2 h^{0,2}_2 -2h^{0,2}_{\bar{\partial}}  
\end{eqnarray*}
We note that these agree with the formula $h^{1,1}_{BC}+h^{2,2}_{BC} = 2h^{1,2}_{BC} + 2$
derived by the author in \cite{McHugh}.
We summarize this with a table:  \newpage
 \begin{center}
 {\bf{Bott-Chern cohomology on complex}} ${\mathbf{S^6}}$ \\
\begin{tabular}{c|c|c|c|c|}
\ $h^{p,q}_{BC}$ & $q= 0$ & $ q=1$ & $q=2$ & $ q=3$ \\
\hline 
$p=0$ &   1   &  0  & $h^{2,0}_{\bar{\partial}} $ & $0$  \\ 
& \  &  \  &  \   & \   \\ 
$p=1$ & 0 & \  \ \  \ \ $2+2h^{0,2}_{\bar{\partial}}$ \ \  & $h^{1,2}_{\bar{\partial}}   + h^{0,2}_2  $ &  $h^{0,2}_{\bar{\partial}} $ \\ 
& \  &  \  &  \   & \   \\ 
$p=2$ & $h^{2,0}_{\bar{\partial}} $ & $ h^{1,2}_{\bar{\partial}}  + h^{0,2}_2 $ &  \ \ $ 2 h^{1,2}_{\bar{\partial}} +2 h^{0,2}_2 -2h^{0,2}_{\bar{\partial}} $ \ \ &\ \ \ $1+h^{0,2}_{\bar{\partial}}+h^{2,0}_{\bar{\partial}}  $ \\ 
& \  &  \  &  \   & \   \\ 
$p=3$ & $0$  &  $h^{0,2}_{\bar{\partial}}  $  &  $1+h^{0,2}_{\bar{\partial}}+h^{2,0}_{\bar{\partial}}  $   & 1  
\end{tabular}
\end{center}
\ \\

\subsection{Bott-Chern Aeppli cohomology of a Calabi-Eckmann 3-fold}
We shall compute the Bott-Chern/Aeppli cohomology for a more concrete case. 
Consider the Calabi-Eckmann 3-fold, diffeomorphic to $S^3 \times S^3$.  We will calculate the Bott-Chern/Aeppli cohomology in this case, replicating the calculation of Angella and Tomassini\cite{AngellaAndTomassini}.  The Dolbeault cohomology of Calabi-Eckmann manifolds was originally calculated by Borel.
For our particular Calabi-Eckmann complex 3-fold ( see for example Cirici\cite{Cirici}) the Dolbeault cohomology is 
\begin{center}
 {\bf{Dolbeault cohomology of the Calabi-Eckmann Complex 3-fold}} \\
\begin{tabular}{c|c|c|c|c|}
\ $h^{p,q}_{\bar{\partial}}$ & $q= 0$ & $ q=1$ & $q=2$ & $ q=3$ \\
\hline 
$p=0$ &   1   &  1  & 0 & 0  \\ 
& \  &  \  &  \   & \   \\ 
$p=1$ & 0 & 1 & 1 &  0 \\ 
& \  &  \  &  \   & \   \\ 
$p=2$ & 0 & 1 &  1  & 0 \\ 
& \  &  \  &  \   & \   \\ 
$p=3$ & 0 &  0  & 1   & 1  
\end{tabular}
\end{center}
Note, that we have $b^1 = b^2 =0$, $h^{2,1}_{\bar{\partial}} =1$ and $h^{0,2}_{\bar{\partial}}=h^{3,1}_{\bar{\partial}} =0$.  Thus $h^{3,1}_2 = 0$.  We also must have $k^{1,2}  =0$ or  $k^{1,2}  =1$.  
We show that $k^{1,2}  =0$.  Indeed, the Dolbeault cohomology in this case (see \cite{AngellaAndTomassini}) has nonzero, $\theta \in H^{0,1}_{\bar{\partial}}$ and nonzero $\partial \theta \in H^{1,1}_{\bar{\partial}}$.   One can check that $H^{2,1}_{\bar{\partial}}$ is generated by ${\bar{\theta}} \partial \theta = \partial ({\bar{\theta}} \theta )$.  Thus the map $H^{2,1}_{\bar{\partial}} \rightarrow H^{2,1}_A$ in our sequence above is the $0$-map and we conclude that $k^{2,1} = k^{1,2} = 0$.
Thus our table for Bott-Chern cohomology is given by 
\begin{center}
 {\bf{Bott-Chern cohomology of the Calabi-Eckmann 3-fold}} \\
\begin{tabular}{c|c|c|c|c|}
\ $h^{p,q}_{BC}$ & $q= 0$ & $ q=1$ & $q=2$ & $ q=3$ \\
\hline 
$p=0$ &   1   &  0  & 0 & 0  \\ 
& \  &  \  &  \   & \   \\ 
$p=1$ & 0 & 2 & 1 &  0 \\ 
& \  &  \  &  \   & \   \\ 
$p=2$ & 0 & 1 &  1  & 1 \\ 
& \  &  \  &  \   & \   \\ 
$p=3$ & 0 &  0  & 1   & 1  
\end{tabular}
\end{center}  
This agrees with the calculation of the Bott-Chern cohomology for the Calabi-Eckmann 3-fold done by Angella and Tomassini\cite{AngellaAndTomassini}.

\subsection{Comparison with Angella's calculation on the Iwasawa manifold and its deformations}
Angella\cite{Angella} gives the De Rham and Dolbeault cohomologies and  calculates the Bott-Chern/Aeppli cohomolgies on the Iwasawa manifold, ${\mathbf I}_3$, and its small deformations in the Kuranishi family of deformations due to Nakamura\cite{Nakamura}.  We will be following this presentation closely, using similiar notation in order to show agreement with our results.  Thus please see Angella\cite{Angella} for fuller details and clarification. 

 The Iwasawa manifold is a holomorphically parallelizable compact complex three-fold with a global holomorphic coframe given by (see Angella\cite{Angella}, p. 39 )
\begin{eqnarray*}
\phi^1  =  dz^1 , \  & \phi^2  =  dz^2 , \  & \phi^3  =  dz^1-z^1dz^2 
\end{eqnarray*}
where $z^1,\,z^2, \, z^3$ are {\em local} coordinates given from ${\mathbf I}_3$ being a quotient space of ${\mathbf{C}}^3$ by a discrete group action. 
The structure equations are easily seen to be
\begin{eqnarray*}
d\phi^1  =  0  , \  & d\phi^2  =  0 , \  & d\phi^3  =  -\phi^1\wedge \phi^2 
\end{eqnarray*}
For small deformations in the Kuranishi family of deformations of ${\mathbf I}_3$, there is a global coframe of $\bigwedge^{1,0}{\mathbf I}_3$ , 
$$\{ \phi^1_t, \phi^2_t,\phi^3_t \} \  , $$ 
where $t=(t_{11}, t_{12}, t_{21}, t_{22}, t_{31}, t_{32}) \in {\bf \Delta}(\epsilon,0) \subset {\mathbf{C}}^6$ and ${\bf \Delta}(\epsilon,0) $ denotes an open ball of radius $\epsilon$ centered at $0$ .  Define 
$D(t) = \begin{pmatrix} 
t_{11} & t_{12} \\
t_{21} & t_{22} 
\end{pmatrix} \ .  $  The different classes of the small Kuranishi deformations of $\mathbf{I}_3$ are defined according to the parameter, $t$: \\
\ \\
$(i): t_{11}= t_{12}= t_{21}= t_{22}=0$, (Note, generically, we have $t_{3k} \neq 0$.) \\
$(ii): ( t_{11}, t_{12}, t_{21}, t_{22}) \neq (0 ,0,0,0), \  D(t) = 0 $ \\
$(iii):  D(t) \neq 0 $ \\
\ \\
Define 
$$ S = 
\begin{pmatrix}
{\overline{\sigma_{1\bar{1}}}} & {\overline{\sigma_{2\bar{2}}}}  & {\overline{\sigma_{2\bar{1}}}}  & {\overline{\sigma_{1\bar{2}}}} \\
 \sigma_{1\bar{1}} &  \sigma_{2\bar{2}} &  \sigma_{2\bar{1}} &  \sigma_{1\bar{2}}
\end{pmatrix}
$$
 where $\sigma_{12}$ and $\sigma_{i{\overline{j}}}$, $i =1,2$, $j=1,2$ will be defined below in the structure equations for classes $ii)$ and $iii)$ (see equation \ref{structeqniiandiii}) \\
The subclasses $(ii.a)$ and $(ii.b)$  of class $(ii)$ are defined by \\
$(ii.a) : D(t) =0$ and $S$ has rank 1. \\
$(ii.b) : D(t) =0$ and $S$ has rank 2. \\
\ \\
Similiarly, the subclasses $(iii.a)$ and $(iii.b)$  of class $(iii)$ are defined by 
\ \\
$(iii.a) : D(t) \neq 0$ and $S$ has rank 1. \\
$(iii.b) : D(t) \neq 0$ and $S$ has rank 2. \\
\ \\
Define $ \zeta^j_t  = z^j+ \sum^2_{k=1} t_{jk} {\bar{z}}^k$ for $j = 1,2$, and 
$$\zeta^3_t = z^3+ \sum^2_{k=1}(t_{3k}+ t_{2k}z^1) {\bar{z}}^k +\frac{1}{2}(t_{11}t_{21}({\bar{z}}^1)^2 +2t_{11}t_{22}{\bar{z}}^1{\bar{z}}^2 +t_{12}t_{22}({\bar{z}}^2)^2) +D(t) {\bar{z}}^3 \ . $$
Thus, define $\phi^j_t = d\zeta^j_t $ for $j=1,2$ and
$$\phi^3_t = d\zeta^3_t - z^1d\zeta^2_t - (t_{21}{\bar{z}}^1 +t_{22}{\bar{z}}^2)d\zeta^1_t \ . $$
The structure equations for class $(i)$ are:
\begin{eqnarray*}
d\phi^1_t  = 0 , \  &  d\phi^2_t  =  0 ,  \  & d\phi^3_t  =  -\phi^1_t \wedge \phi^2_t 
\end{eqnarray*}
For classes $(ii)$ and $(iii)$, the structure equations are given by
\begin{eqnarray} \label{structeqniiandiii}
d\phi^1_t  & = & 0 , \ \ \ \  d\phi^2_t  =  0 ,  
\end{eqnarray}
\begin{eqnarray*}
d\phi^3_t  & = & \sigma_{12} \phi^1_t\wedge \phi^2_t   + \sigma_{1\bar{1}} \phi^1_t\wedge {\bar{\phi^1_t }}+ \sigma_{1 \bar{2}} \phi^1_t\wedge {\bar{\phi^2_t} } 
 \\
  \ & \ & + \sigma_{2\bar{1}} \phi^2_t\wedge {\bar{\phi^1_t}} + \sigma_{2\bar{2}} \phi^2_t\wedge {\bar{\phi^2_t}} 
\end{eqnarray*}
where $\sigma_{12} \in {\mathbf{C}}$ and $\sigma_{i{\overline{j}}} \in {\mathbf{C}}$, $i =1,2$, $j=1,2$ are complex numbers depending only on $t$.

We shall be using notation as in Angella\cite{Angella}: $ \phi^{I{\overline{J}}} = \phi^{i_1} \wedge \cdots \wedge \phi^{i_p}\wedge  {\overline{\phi^{j_1}}} \wedge \cdots \wedge {\overline{\phi^{j_q}}}$ for multi-indices, $I$ and $J$. \\
\ \\
One can check that $h^{3,1}_2=2$ for all the classes of deformations of ${\mathbf{I}}_3$:  One can check that, a priori, $H^{3,1}_{\bar{\partial}}= ker(\partial )$ and that, $$H^{3,1}_{\bar{\partial}} = {\mathbf{C}}< \phi_t^{123{\bar{1}}},\phi_t^{123{\bar{2}}}>  \ .$$  Since, $H^{2,1}_{\bar{\partial}} \subseteq {\mathbf{C}}< \phi_t^{12{\bar{1}}},\phi_t^{12{\bar{2}}},\phi_t^{13{\bar{1}}},\phi_t^{13{\bar{2}}},\phi_t^{23{\bar{1}}},\phi_t^{23{\bar{2}}}>  $ for class $(i)$ and 
$$H^{2,1}_{\bar{\partial}} \subseteq {\mathbf{C}}< \phi_t^{12{\bar{1}}},\phi_t^{12{\bar{2}}},\phi_t^{13{\bar{1}}}-\frac{\sigma_{2{\bar{2}}}}{{\overline{\sigma_{12}}}}\phi_t^{12{\bar{3}}},\phi_t^{13{\bar{2}}}-\frac{\sigma_{2{\bar{1}}}}{{\overline{\sigma_{12}}}}\phi_t^{12{\bar{3}}},  \ \ \ \ \ \ \ \ \ \ \ \ \ \ \ \  \ \ \ \ \ \  \ \ \ \ \ \ \ $$
$$ \ \ \ \ \ \ \ \ \ \ \ \ \ \ \ \ \ \ \ \ \ \ \ \ \ \ \ \ \ \ \ \ \ \ \ \ \ \ \ \ \ \ \ \phi_t^{23{\bar{1}}}-\frac{\sigma_{1{\bar{2}}}}{{\overline{\sigma_{12}}}}\phi_t^{12{\bar{3}}},\phi_t^{23{\bar{2}}}-\frac{\sigma_{1{\bar{1}}}}{{\overline{\sigma_{12}}}}\phi_t^{12{\bar{3}}}>  $$
for classes $(ii)$ and $(iii)$,
we see that $\partial$ acts on $H^{2,1}_{\bar{\partial}}$ as the zero map.  Thus $Im(\partial : H^{2,1}_{\bar{\partial}} \rightarrow H^{3,1}_{\bar{\partial}})  =0$ and 
$h^{3,1}_2=2$.

We see from the $2,1$ -forms above, that $\phi_t^{12{\bar{1}}} =\partial \phi_t^{3,{\bar{1}}}$ and $\phi_t^{12{\bar{2}}} =\partial \phi_t^{3,{\bar{2}}}$ and thus $k^{2,1}\leq 4$. We shall show  that
$k^{1,2}=4$ on the Iwasawa manifold and its classes of deformations {\em (i), (ii)} and {\em (iii)}.  In fact, Angella\cite{Angella} shows that the  $2,1$ -forms, 
$$\{ \phi_t^{13{\bar{1}}},\phi_t^{13{\bar{2}}},\phi_t^{23{\bar{1}}},\phi_t^{23{\bar{2}}} \}  $$ for class $(i)$ and 
$$\{ \phi_t^{13{\bar{1}}}-\frac{\sigma_{2{\bar{2}}}}{{\overline{\sigma_{12}}}}\phi_t^{12{\bar{3}}},\phi_t^{13{\bar{2}}}-\frac{\sigma_{2{\bar{1}}}}{{\overline{\sigma_{12}}}}\phi_t^{12{\bar{3}}}, \phi_t^{23{\bar{1}}}-\frac{\sigma_{1{\bar{2}}}}{{\overline{\sigma_{12}}}}\phi_t^{12{\bar{3}}},\phi_t^{23{\bar{2}}}-\frac{\sigma_{1{\bar{1}}}}{{\overline{\sigma_{12}}}}\phi_t^{12{\bar{3}}} \}  $$
for classes $(ii)$ and $(iii)$, are ${\bar{\partial}}_{J_t}$-harmonic with respect to the hermitian metric $g_t = \phi_t^1 \odot{\bar{ \phi_t^1}} +\phi_t^2 \odot {\bar{\phi_t^2}}+\phi_t^3 \odot {\bar{\phi_t^3}}$.

We can show that for the class of deformations $(i)$, the $2,1$-forms
$$\{ \phi_t^{13{\bar{1}}},\phi_t^{13{\bar{2}}},\phi_t^{23{\bar{1}}},\phi_t^{23{\bar{2}}} \}  $$ are also $\partial_{J_t}$-harmonic:
It is easy to see from the structure equations that these are $\partial$-closed.  Now observe that  (using the notation of Wells\cite{Wells})
\begin{eqnarray*}
{\bar{\ast}}  \phi_t^{13{\bar{1}}}  =  \phi_t^{2{\bar{2}\bar{3}}} ,  \  & {\bar{\ast}}  \phi_t^{13{\bar{2}}}  =  - \phi_t^{2{\bar{1}\bar{3}}} \\
{\bar{\ast}}  \phi_t^{23{\bar{1}}}  =  - \phi_t^{1{\bar{2}\bar{3}}} , \  & {\bar{\ast}}  \phi_t^{23{\bar{2}}}  =   \phi_t^{1{\bar{1}\bar{3}}} \\
\end{eqnarray*}
and that these are also $\partial$-closed.  Thus these $2,1$-forms for class $(i)$- are $\partial^{*}$-closed and hence $\partial_{J_t}$-harmonic.  Since they are also  ${\bar{\partial}}_{J_t}$-harmonic, we know that  ${\mathbf{C}} < \phi_t^{13{\bar{1}}},\phi_t^{13{\bar{2}}},\phi_t^{23{\bar{1}}},\phi_t^{23{\bar{2}}} > $ is orthogonal to $Im({\bar{\partial}}) + Im(\partial )$.   We conclude that $k^{2,1} \geq 4$ for the class $(i)$.  Thus, $k^{2,1} = 4$ for class $(i)$.

For deformation classes, $(ii)$ and $(iii)$, consider the five dimensional complex vector space,
$$W={\mathbf{C}} < \phi_t^{12{\bar{3}}} , \phi_t^{13{\bar{1}}},\phi_t^{13{\bar{2}}},\phi_t^{23{\bar{1}}},\phi_t^{23{\bar{2}}} >  \ . $$
One can check that this is also
$$W={\mathbf{C}} <{\bar{ \partial}}^*\phi_t^{12{\bar{1}}{\bar{2}}}, \phi_t^{13{\bar{1}}}-\frac{\sigma_{2{\bar{2}}}}{{\overline{\sigma_{12}}}}\phi_t^{12{\bar{3}}},\phi_t^{13{\bar{2}}}-\frac{\sigma_{2{\bar{1}}}}{{\overline{\sigma_{12}}}}\phi_t^{12{\bar{3}}}, \ \ \ \ \ \ \ \ \ \ \ \ \ \ \ \  \ \ \ \ \ \  \ \ \ \ \ \ \ $$
$$ \ \ \ \ \ \ \ \ \ \ \ \ \ \ \ \ \ \ \ \ \ \ \ \ \ \ \ \ \ \ \ \ \ \ \ \ \ \ \ \ \ \ \ \phi_t^{23{\bar{1}}}-\frac{\sigma_{1{\bar{2}}}}{{\overline{\sigma_{12}}}}\phi_t^{12{\bar{3}}},\phi_t^{23{\bar{2}}}-\frac{\sigma_{1{\bar{1}}}}{{\overline{\sigma_{12}}}}\phi_t^{12{\bar{3}}}>  \ . $$
Recalling the orthogonal decomposition, $ {\mathcal{E}}^{2,1} = {\mathcal{H}}^{2,1}_{\bar{\partial}} \oplus Im({\bar{\partial}}) \oplus  Im({\bar{\partial}}^* )  \ , $
we see that $W$ is orthogonal to $Im({\bar{\partial}})$.

We also have
$$W={\mathbf{C}} < \partial \phi_t^{3{\bar{3}}}, \phi_t^{13{\bar{1}}}+\frac{\overline{\sigma_{1{\bar{1}}}}}{\sigma_{12}}\phi_t^{12{\bar{3}}},\phi_t^{13{\bar{2}}}+\frac{\overline{\sigma_{2{\bar{1}}}}}{\sigma_{12}}\phi_t^{12{\bar{3}}},  \ \ \ \ \ \ \ \ \ \ \ \ \ \ \ \  \ \ \ \ \ \  \ \ \ \ \ \ \ $$
$$ \ \ \ \ \ \ \ \ \ \ \ \ \ \ \ \ \ \ \ \ \ \ \ \ \ \ \ \ \ \ \ \ \ \ \ \ \ \ \ \ \ \ \ \phi_t^{23{\bar{1}}}+\frac{\overline{\sigma_{1{\bar{2}}}}}{\sigma_{12}}\phi_t^{12{\bar{3}}},\phi_t^{23{\bar{2}}}+\frac{\overline{\sigma_{2{\bar{2}}}}}{\sigma_{12}}\phi_t^{12{\bar{3}}}>  $$
noting that for the first basis element we have
 $$ \partial \phi_t^{3{\bar{3}}} = \sigma_{12} \phi_t^{12{\bar{3}}} + {\overline{\sigma_{1{\bar{1}}}}}\phi_t^{13{\bar{1}}}+ {\overline{\sigma_{2{\bar{1}}}}}\phi_t^{13{\bar{2}}} + {\overline{\sigma_{1{\bar{2}}}}}\phi_t^{23{\bar{1}}}+ {\overline{\sigma_{2{\bar{2}}}}}\phi_t^{23{\bar{2}}}$$
and that the other four basis elements are $\partial$-harmonic.
Thus we see that $W \cap Im(\partial ) = {\mathbf{C}}< \partial \phi_t^{3{\bar{3}}}>$.  
We also have that 
$$ W = {\mathbf{C}} <\frac{1}{\sigma_{12}} \partial \phi_t^{3{\bar{3}}},  \phi_t^{13{\bar{1}}}-\frac{\sigma_{2{\bar{2}}}}{{\overline{\sigma_{12}}}}\phi_t^{12{\bar{3}}},\phi_t^{13{\bar{2}}}-\frac{\sigma_{2{\bar{1}}}}{{\overline{\sigma_{12}}}}\phi_t^{12{\bar{3}}},  \ \ \ \ \ \ \ \ \ \ \ \ \ \ \ \  \ \ \ \ \ \  \ \ \ \ \ \ \ $$
$$ \ \ \ \ \ \ \ \ \ \ \ \ \ \ \ \ \ \ \ \ \ \ \ \ \ \ \ \ \ \ \ \ \ \ \ \ \ \ \ \ \ \ \ \phi_t^{23{\bar{1}}}-\frac{\sigma_{1{\bar{2}}}}{{\overline{\sigma_{12}}}}\phi_t^{12{\bar{3}}},\phi_t^{23{\bar{2}}}-\frac{\sigma_{1{\bar{1}}}}{{\overline{\sigma_{12}}}}\phi_t^{12{\bar{3}}}>   $$
since for small $t$, the matrix,
$$ \begin{pmatrix}
1 &  -\frac{\sigma_{2{\bar{2}}}}{\overline{\sigma_{12}}} &  -\frac{\sigma_{2{\bar{21}}}}{\overline{\sigma_{12}}} & -\frac{\sigma_{1{\bar{2}}}}{\overline{\sigma_{12}}} &  -\frac{\sigma_{1{\bar{1}}}}{\overline{\sigma_{12}}} \\ 
\frac{\overline{\sigma_{1{\bar{1}}}}}{\sigma_{12}} & 1 & 0 & 0 & 0 \\ \\
\frac{\overline{\sigma_{2{\bar{1}}}}}{\sigma_{12}} & 0 & 1 & 0 & 0 \\ \\
\frac{\overline{\sigma_{1{\bar{2}}}}}{\sigma_{12}} & 0 & 0 & 1 & 0 \\ \\
\frac{\overline{\sigma_{2{\bar{2}}}}}{\sigma_{12}} & 0 & 0 & 0 & 1 
\end{pmatrix} $$
has non-zero determinant.  Note that the basis elements 
$$ \{ \phi_t^{13{\bar{1}}}-\frac{\sigma_{2{\bar{2}}}}{{\overline{\sigma_{12}}}}\phi_t^{12{\bar{3}}},\phi_t^{13{\bar{2}}}-\frac{\sigma_{2{\bar{1}}}}{{\overline{\sigma_{12}}}}\phi_t^{12{\bar{3}}}, \phi_t^{23{\bar{1}}}-\frac{\sigma_{1{\bar{2}}}}{{\overline{\sigma_{12}}}}\phi_t^{12{\bar{3}}},\phi_t^{23{\bar{2}}}-\frac{\sigma_{1{\bar{1}}}}{{\overline{\sigma_{12}}}}\phi_t^{12{\bar{3}}} \}   $$
are ${\bar{\partial}}$-closed and $\partial$-closed.  Denote,
$$V ={\mathbf{C}}< \phi_t^{13{\bar{1}}}-\frac{\sigma_{2{\bar{2}}}}{{\overline{\sigma_{12}}}}\phi_t^{12{\bar{3}}},\phi_t^{13{\bar{2}}}-\frac{\sigma_{2{\bar{1}}}}{{\overline{\sigma_{12}}}}\phi_t^{12{\bar{3}}}, \phi_t^{23{\bar{1}}}-\frac{\sigma_{1{\bar{2}}}}{{\overline{\sigma_{12}}}}\phi_t^{12{\bar{3}}},\phi_t^{23{\bar{2}}}-\frac{\sigma_{1{\bar{1}}}}{{\overline{\sigma_{12}}}}\phi_t^{12{\bar{3}}} >$$
Thus, $V / (V \cap (Im({\bar{\partial}}) + Im(\partial ) )) \subseteq K^{2,1} \  . $
Since 
$$V / (V \cap (Im({\bar{\partial}}) + Im(\partial ) ))  = V / (V \cap {\mathbf{C}}<\partial \phi^{3{\bar{3}}}> )$$
we have that $dim_{\mathbf{C}}(V / (V \cap (Im({\bar{\partial}}) + Im(\partial ) )) )=4$ and $k^{2,1} \geq 4$.  We conclude that $k^{2,1} =4$ for deformation classes, $(ii)$ and $(iii)$.

It is clear from $H^{0,1}_{\bar{\partial}} = {\mathbf{C}}<{\bar{\phi_t}}^1, {\bar{\phi_t}}^2 >$
 that $k^{1,0}=k^{0,1}=2$. 

We now calculate $k^{2,0}$ for each of the three classes.  For class $(i)$, we have that 
$$H^{2,0}_{\bar{\partial}} = {\mathbf{C}}<\phi_t^{12},\phi_t^{13},\phi_t^{23}>$$
where we note that these basis elements are also $\partial$-closed and that $\phi_t^{12} = \partial \phi_t^3$.  Thus for class $(i)$, $k^{2,0}=2$.  

For class $(ii)$, we have that 
$$H^{2,0}_{\bar{\partial}} = {\mathbf{C}}<\phi_t^{12},\alpha \phi_t^{13}+\beta \phi_t^{23}>$$
with $\alpha$ and $\beta$ not both zero.  As before, the generators are also $\partial$-closed and  $\phi_t^{12} = \partial \phi_t^3$.  Thus for class $(ii)$, $k^{2,0}=1$.

For class $(iii)$, we have that 
$H^{2,0}_{\bar{\partial}} = {\mathbf{C}}<\phi_t^{12}>$
and as before,  $\phi_t^{12} = \partial \phi_t^3$.  Thus for class $(iii)$, $k^{2,0}=0$.\\

We now calculate $k^{1,1}$ for each of the three classes. For class $(i)$, $K^{1,1}$ has a basis consisting of $\{ \phi_t^{1\overline{1}},  \phi_t^{1\overline{2}},  \phi_t^{2\overline{1}},  \phi_t^{2\overline{2}}\}$. Thus, $k^{1,1}  = 4$ for class $(i)$.    

For classes $(ii. a)$ and $(iii. a)$, $${\bar{\partial}} \phi_t^3 =  \sigma_{1\bar{1}} \phi_t^{1\bar{1}}+ \sigma_{1 \bar{2}} \phi^{1\bar{2}}_t + \sigma_{2\bar{1}} \phi^{2\bar{1}}_t + \sigma_{2\bar{2}} \phi^{2\bar{2}}_t $$
and 
$$ \partial {\overline{ \phi_t^3}} =-( \overline{ \sigma_{1\bar{1}}} \phi_t^{1\bar{1}}+\overline{ \sigma_{2 \bar{1}}} \phi^{1\bar{2}}_t + \overline{\sigma_{1\bar{2}}} \phi^{2\bar{1}}_t +\overline{ \sigma_{2\bar{2}}} \phi^{2\bar{2}}_t) $$
are non-zero but linearly dependent over $\mathbf{C}$ (since $S$ has rank 1) and thus a non-zero linear combination of  $\{ \phi_t^{1\overline{1}},  \phi_t^{1\overline{2}},  \phi_t^{2\overline{1}},  \phi_t^{2\overline{2}}\}$ is in $Im(\partial) + Im({\bar{\partial}})$.  Hence,  $k^{1,1}=3$ in classes $(ii. a)$ and $(iii. a)$.  

For classes $(ii. b)$ and $(iii. b)$, $${\bar{\partial}} \phi_t^3 =  \sigma_{1\bar{1}} \phi_t^{1\bar{1}}+ \sigma_{1 \bar{2}} \phi^{1\bar{2}}_t + \sigma_{2\bar{1}} \phi^{2\bar{1}}_t + \sigma_{2\bar{2}} \phi^{2\bar{2}}_t $$
and 
$$ \partial{\overline{ \phi_t^3}} =-( \overline{ \sigma_{1\bar{1}}} \phi_t^{1\bar{1}}+\overline{ \sigma_{2 \bar{1}}} \phi^{1\bar{2}}_t + \overline{\sigma_{1\bar{2}}} \phi^{2\bar{1}}_t +\overline{ \sigma_{2\bar{2}}} \phi^{2\bar{2}}_t )$$
are non-zero and linearly independent over $\mathbf{C}$.  Thus two linearly independent non-zero linear combinations of  $\{ \phi_t^{1\overline{1}},  \phi_t^{1\overline{2}},  \phi_t^{2\overline{1}},  \phi_t^{2\overline{2}}\}$ are in $Im(\partial) + Im({\bar{\partial}})$.  Hence, $k^{1,1}=4-2=2$ in classes $(ii. b)$ and $(iii. b)$. 

Clearly, $k^{3,0} = h^{0,3}_{\bar{\partial}} = 1$ for all the classes of deformations of the Iwasawa manifold being considered. 
We summarize the $k^{i,j}$ with the table,  
\begin{center}  {\bf{Table of $k^{i,j}$ for Classes of Deformations of the Iwasawa manifold}} \\
  \begin{tabular}{ c | c c c c c }
    \hline
    $k^{i,j}$ & $k^{1,0}$ & $k^{1,1}$ &$k^{1,2}$ & $k^{2,0}$  & $k^{3,0}$  \\ \hline
  $(i)$ & 2 & 4 & 4 & 2 & 1 \\ \\
   $(ii.a)$ & 2 & 3 & 4 & 1 & 1   \\ \\
 $(ii.b)$ & 2 & 2 & 4 & 1& 1  \\ \\
 $(iii.a)$ & 2 & 3 & 4 & 0 & 1   \\ \\
 $(iii.b)$ & 2 & 2 & 4 & 0 & 1   \\
    \hline
  \end{tabular}
\end{center}
We also include here Angella's calculation of the Dolbeault cohomology for the classes of deformations of the Iwasawa manifold, \\
\begin{center}  {\bf{Table of $h_{\bar{\partial}}^{i,j}$ for Classes of Deformations of the Iwasawa manifold}} \\
  \begin{tabular}{ c | c c c c c c c }
    \hline
    $h_{\bar{\partial}}^{i,j}$ & $h_{\bar{\partial}}^{1,0}$ & $h_{\bar{\partial}}^{0,1}$ & $h_{\bar{\partial}}^{1,1}$ &$h_{\bar{\partial}}^{1,2}$ & $h_{\bar{\partial}}^{0,2}$ & $h_{\bar{\partial}}^{2,0}$ & $h_{\bar{\partial}}^{3,0}$  \\ \hline
    $(i)$ & 3 & 2 & 6 & 6 & 2 & 3 & 1 \\ \\
 $(ii.a)$ & 2 & 2 & 5 & 5 & 2 & 2 & 1 \\ \\
 $(ii.b)$ & 2 & 2 & 5 & 5 & 2 & 2 & 1 \\ \\
 $(iii.a)$ & 2 & 2 & 5 & 4 & 2 & 1 & 1 \\ \\
 $(iii.b)$ & 2 & 2 & 5 & 4 & 2 & 1 & 1 \\
    \hline
  \end{tabular}
\end{center}

Let us calculate $h^{1,1}_{BC}$ for all the classes of small Kuranishi deformations of ${\mathcal{I}}_3$.  We have
$$h^{1,1}_{BC} =   2h^{0,1}_{\bar{\partial}} - b^1 +dim_{\mathbf{R}}(H^{1,1}_{dR}({\mathbf{R}})) \ . $$ 
One can check that for all five classes of deformations we have that 
$$ {\mathbf{R}} < i\phi^{1{\bar{1}}}_t,  i\phi^{2{\bar{2}}}_t,  \phi^{1{\bar{2}}}_t- \phi^{2{\bar{1}}}_t, i(\phi^{1{\bar{2}}}_t+ \phi^{2{\bar{1}}}_t) > \ 
\  \subseteq \  H^{1,1}_{dR}({\mathbf{R}}) \  . $$
Thus, $dim_{\mathbf{R}}( H^{1,1}_{dR}({\mathbf{R}}) ) \geq 4 $.   We wish to show $dim_{\mathbf{R}}( H^{1,1}_{dR}({\mathbf{R}}) ) = 4 $ for all the classes of deformations.  
 Any $d$-closed $1,1$-form must also be ${\bar{\partial}}$-closed and $\partial$-closed.   Such a $1,1$-form independent of 
$$ {\mathbf{R}} < i\phi^{1{\bar{1}}}_t,  i\phi^{2{\bar{2}}}_t,  \phi^{1{\bar{2}}}_t- \phi^{2{\bar{1}}}_t, i(\phi^{1{\bar{2}}}_t+ \phi^{2{\bar{1}}}_t) > $$
will be written as 
$$ \psi = \alpha \phi^{1\bar{3}}+\beta \phi^{2\bar{3}}+\gamma \phi^{3\bar{1}}+ \delta \phi^{3\bar{2}}  + \epsilon \phi^{3\bar{3}} $$
where the $\alpha , \beta , \gamma , \delta $ and $\epsilon$ are global smooth complex valued functions.  Applying ${\bar{\partial}}$ and $\partial$-closedness and projecting onto 
$${\mathbf{C}}<\phi^{1{\bar{1}}{\bar{2}}},\phi^{2{\bar{1}}{\bar{2}}},\phi^{12{\bar{1}}}, \phi^{12{\bar{2}}}>$$
gives the matrix equation \\
$$ \begin{pmatrix}
-{\overline{\sigma_{12}}} & 0 & - \sigma_{1{\bar{2}}} &  - \sigma_{1{\bar{1}}} \\
0 & -{\overline{\sigma_{12}}}  & - \sigma_{2{\bar{2}}} &  - \sigma_{2{\bar{1}}} \\
  {\overline{ \sigma_{1{\bar{2}}}}} &  -{\overline{ \sigma_{1{\bar{1}}}}}& \sigma_{12} & 0  \\
 {\overline{\sigma_{2{\bar{2}}}}} &  -{\overline{ \sigma_{2{\bar{1}}}}}& 0 & \sigma_{12}  
\end{pmatrix} 
\begin{pmatrix}
\alpha \\
\beta \\
\gamma \\
\delta 
\end{pmatrix} = 
\begin{pmatrix}
0 \\
0 \\
0 \\
0
\end{pmatrix} \ \ . $$
We note that this is similiar to the equation in Angella\cite{Angella}, p.47.   The matrix above has rank 4 for small deformations in classes $(i)$, $(ii)$, and $(iii)$.   Thus we must have $\alpha = \beta = \gamma = \delta = 0$.  We then can easily show that  $\epsilon = 0$:
$$ \partial (\epsilon \phi^{3 {\bar{3}}} ) = (\partial \epsilon ) \wedge \phi^{3 {\bar{3}}} - \epsilon \phi^{12{\bar{3}}} + \epsilon \phi^3 \wedge (\partial {\overline{\phi^3}}) = 0  \  . $$
Each of these three terms are linearly independent, i.e. the only term in $ \phi^{12{\bar{3}}}$, when all expanded,  is the middle term.  Thus we must have
$\epsilon = 0$.   Hence we conclude that for the Iwasawa manifold and its small deformations in the Kuranishi family,
$dim_{\mathbf{R}}( H^{1,1}_{dR}({\mathbf{R}}) ) = 4 $.  One can see that we obtain for all the classes of deformation $(i)$, $(ii)$, and $(iii)$,
\begin{eqnarray*}
h^{1,1}_{BC}& =& 2h^{0,1}_{\bar{\partial}} - b^1 +dim_{\mathbf{R}}( H^{1,1}_{dR}({\mathbf{R}}) )  =  2*2-4+4  =  4
\end{eqnarray*}
This agrees with Angella\cite{Angella}p.49.

Using this value of $h^{1,1}_{BC}$, the values for $k^{i,j}$ in the table above,  Angella's calculations of  $h^{p,q}_{\bar{\partial}}$ in the above table for $h^{i,j}_{BC}$, we get agreement with Angella's calculation of the rest of the Bott-Chern cohomology.
For example,
\begin{eqnarray*}
h^{2,2}_{BC} &=& -  h^{1,1}_{BC} +h^{0,1}_{\bar{\partial}} -h^{0,2}_{\bar{\partial}} -h^{1,0}_{\bar{\partial}}+ h^{1,1}_{\bar{\partial}}+ h^{1,2}_{\bar{\partial}}  +h^{2,0}_{\bar{\partial}} \\
  \ & \  & \ \ \ \ \ \ \ \ \ \ \ \ \ \ \ \ \ \ \ \ \ \ \ \ \ \ \ \ \ \ \ \ \ \ \ \ \ \ \ +2 h^{3,1}_2 +k^{1,1}-k^{1,2}-2k^{2,0} \ .
\end{eqnarray*}
For the different classes of deformations, the reader can check that this is:
\begin{center} 
  \begin{tabular}{ c | c c c c c }
    \hline
Class &\ \ \ $(i)$ \ \ \ \ &  $(iia) \  \ \  $ &   $(iib)$  \ \ \  &  $(iiia)$ \ \ \   & $(iiib)$ \ \ \  \\ 
  $h^{2,2}_{BC}$ & \  8 \  & 7 \ & 6\  &  7 \  & 6 \\ 
    \hline
  \end{tabular}
\end{center}

 The reader can also check the calculation of $h^{1,2}_{BC}$, for each of the three classes $(i)$, $(ii)$, and $(iii)$ using the equation,
$h^{1,2}_{BC}  =   h^{1,2}_{\bar{\partial}}  + h^{3,1}_2 - k^{2,0} $ which produces the following values:

\begin{center} 
  \begin{tabular}{ c | c c c  }
    \hline
Class &\ \ \ $(i)$ \ \ \ \ &  $(ii)$   \ \ \  &  $(iii)$ \ \ \   \\ 
  $h^{1,2}_{BC}$ & \  6 \  & 6 \ & 6 \   \\ 
    \hline
  \end{tabular}
\end{center}

\end{document}